\documentclass[12pt]{article}

\usepackage[utf8]{inputenc}
\usepackage[tbtags]{amsmath}
\usepackage{amsfonts,amssymb,amsthm}
\usepackage{color,graphics}

\usepackage{graphicx}

\numberwithin{equation}{section}

\newtheorem{counter_main}{Empty counter}[section]
\newtheorem{theorem}[counter_main]{Theorem}
\newtheorem{proposition}[counter_main]{Proposition}
\newtheorem{exam}[counter_main]{Example}
\newtheorem{cons}[counter_main]{Corollary}
\newtheorem{lemma}[counter_main]{Lemma}
\newtheorem{problem}[counter_main]{Problem}
\newtheorem{note}[counter_main]{Remark}
\newtheorem{hyp}[counter_main]{Conjecture}
\newtheorem{denote}[counter_main]{Notation}
\newtheorem{defn}[counter_main]{Definition}

\newenvironment{prf}{{\scshape Proof.}}{$\Box$}

\newcommand*{\per}{{\rm{per}}\;}
\newcommand*{\mper}{{\rm{mper}}\;}

\newcommand{\suml}[2]{\sum\limits_{{#1}}^{{#2}}}                                        

\newcommand{\rk}{{\rm{rk}}\,}
\newcommand{\ol}{\overline}

\newcommand{\F}{\mathbb{F}}

\newcommand{\N}{\mathbb{N}}

\newcommand{\al}{\alpha}
\newcommand{\be}{\beta}

\newcommand{\lep}{\le^*}
\newcommand{\gep}{\ge^*}
\def\SS{\Lambda}

\begin{document}

\title{  Kr\"auter conjecture on permanents is true\footnote{The work is  financially supported by the grant  RSF 17-11-01124 }}
\author{M.V. Budrevich, A.E. Guterman}
\date{ {\small Moscow State University,  Moscow, GSP-1, 119991, 
Russia \\ Moscow Institute of Physics and Technology, Dolgoprudny, 141701, Russia }}
\maketitle
\begin{abstract}

In this paper we   investigate the permanent of $(-1,1)$-matrices over fields of zero characteristics and our main goal is to provide a sharp upper bound for the value of the permanent of such matrices depending on matrix  rank, solving  Wang's problem  posed in 1974  by confirming Kr\"auter conjecture formulated in 1985.

\end{abstract}

\section{Introduction}

The main object of our consideration is the permanent function.
\begin{defn}
\label{dfper}
For a square matrix $A$ of size $n$, the permanent function is 
$\per A=\sum_{\sigma \in S_n} a_{1\sigma(1)}\cdots  a_{n\sigma(n)},$
where $S_n$ denotes the permutation group on the set $\{1,2,\ldots,n\}$.
\end{defn}

Permanent is very useful in discrete mathematics, combinatorics and linear algebra, however, if we compare it with the determinant, it is quite difficult to compute this function. In fact, the computation of permanent is a known $\sharp P$-complete problem, see \cite{Minc,Valiant} for the detailed and self-contained information. There are many open problems related to the permanent function. In particular, many problems are related to the behavior of the permanent function on  matrices with specific entries. The central role in this investigations is played by  $(0,1)$- and $(-1,1)$-matrices, which are important for a number of pure and applied problems.

The detailed and self-contained information about the behavior of permanent on  $(0,1)$-matrices  and their applications can be found for example  in \cite{BR}. 
$(-1,1)$-matrices also arise naturally in different economical problems, see \cite{BapatR}, 
in different problems of sign conversion, see \cite{McCr},  
  and in the theory of Hadamard matrices, see \cite{Minc-update}. 

Observe that the investigation of permanents of $(-1,1)$-matrices are much more subtle than  for permanents of $(0,1)$-matrices and other nonnegative matrices. Many standard properties of $(0,1)$-matrices are not true for $(-1,1)$-matrices. For example, the permanent of a $(0,1)$-matrix cannot decrease, if a zero is replaced by one, which is not true if $-1$ is replaced by 1 in $(-1,1)$-matrix. Moreover, the permanent of a $(-1,1)$-matrix is totally independent of the number of $-1$s in a matrix. In fact, it is easy to see that a $(-1,1)$-matrix of even order can attain its maximum permanent value $n!$ even if all the entries are $-1$ or exactly half of the entries are  $-1$, see \cite[Example 1]{Wang}. 

The first investigations of permanents on $(-1,1)$-matrices were started in  \cite{MN, Reich}, later they were continued in the papers \cite{Kr,Perfect,Wang}, see also their bibliography. 
For   results   on  $(-1,1)$-matrices and their permanents we recommend the survey paper~\cite{Wanless} and references therein. More recent develop\-ments and open problems involving permanent can be found in \cite{Bapat,Minc-update,Zhang}. 

Since it is difficult to compute permanent, different estimates for this function   play a central role in permanent theory. As an example, one can consider  van der Waerden conjecture claiming the bound for permanent on the set of doubly stochastic matrices, \cite{Minc}, which stimulated a lot of research on permanents throughout the century before the conjecture was confirmed.

In this paper we   investigate the permanent of $(-1,1)$-matrices over fields of zero characteristics and our main goal is to provide a sharp upper bound for the value of the permanent of such matrices depending on matrix  rank. This result solves the open problem posed by    Wang in 1974, see  \cite{Wang},  by confirming Kr\"auter conjecture (1985) stated in \cite{Kr}, see the exact formulations below. 

Let $\F$ be a field of zero characteristics. We denote by $M_{k,n}(\F)$ the set of $k$-by-$n$ matrices with the coefficients from $\F$.  $M_{k,n}(\pm 1)$  denotes its subset consisting of all matrices whose entries  are $\pm 1$ only. As usual, $M_n=M_{n,n}$ and $GL_n$ denotes the subset  of all  invertible matrices in $M_n$.

The following   problem related to the permanent of $(-1,1)$-matrices was posed in 1974.

\begin{problem} (Wang, {\cite[Problem 2]{Wang}})
Let $A \in M_n(\pm 1)$ be a nonsingular matrix. Is there a sharp upper bound for $|\per A|$?
\label{p1}
\end{problem}

For the detailed information on this problem see also the survey paper by Minc, \cite[Problem 6, p. 245]{Minc0}. Several  interesting  upper bounds were obtained in original paper by Wang and in subsequent papers.  For example, let us present here the result by Wang:

\begin{theorem} (Wang, {\cite[Propositions 2, 3, p. 355]{Wang}}) Let $A \in M_n(\pm 1)$. Then  $\per A\le n!$. Moreover,  $\per A = n!$  if and only if  $A$ can be reduced to $J_n$  (the matrix of all 1)  by    transposition,  row or column permutations and    multiplication of rows or columns  by  $ (- 1) $. Otherwise  $\per A\le  (n-2) (n-1)!$.
\label{th-kr-new}\end{theorem}

However, this bound is not sharp and is not depending on rank. Minc in  \cite{Minkb} presented an upper bound which is also not sharp, see \cite[Formula (22)]{Minc0}. Several  bounds for matrices of special structure (limited numbers  of $-1$s in rows and columns) were obtained by Kr\"auter and Seifter in \cite{KrSe2,KrSe,Seifter} by means  of the permanents of  the following special matrices.
\begin{denote}
Let $D_{(n,k,l)}$ denote the matrix  $D_{(n,k,l)} = (d_{ij}) \in M_{k,n}(\pm 1)$ where $0 \le l \le k \le n$ and
$$
d_{ij} = 
\begin{cases}
-1, i = j \text{ and } j \in \{ 1, \ldots, l\} \\
1, \text{ otherwise.}
\end{cases}
$$
If $n = k$ we write $D_{(n,n,l)} = D_{(n,l)}$.
\end{denote}
 A detailed exposition of current achievements, including methods and difficulties is given in \cite[Section 6]{Kr2}.

Kr\"auter in the paper  \cite{Kr} formulated  the following conjecture which answers the Wang's question and  provides a possible sharp upper bound for the values of the permanent function for matrices from  $M_n(\pm 1)$  via the rank function. Kr\"auter conjecture also appears in Minc's collections of problems related to permanent, see \cite[Conjecture 37,  p. 134]{Minc1}.

\begin{hyp} (Kr\"auter, \cite[Conjecture 5.2, p. 13]{Kr},  \cite[Conjeture 6.2, p. 74]{Kr2}, Minc, \cite[Conjecture 37]{Minc1})
Let $A \in M_n(\pm 1)$ where  $n \ge 5$ and  $\rk A = r+1$ for some $r$, $0\le r\le n-1$. Then the following inequality holds
\begin{equation}
|\per A| \le \per D_{(n,r)}.
\label{eq-27.16}\end{equation}
The equality in 
\eqref{eq-27.16} holds if and only if the matrix  $A$ can be obtained from  $D_{(n,r)}$ by the transposition,  row or column permutations and  multiplications of rows or columns by $ (- 1) $. 
\label{h1/1}
\end{hyp}

For the convenience of notations in the inductive proof of Conjecture  \ref{h1/1} we extend its formulation in the following way.

\begin{hyp}
Let $A \in M_n(\pm 1)$, $n\ge 2$, and  $\rk A = r+1$ for some $r$, $0\le r\le n-1$. Then the inequality 
\begin{equation}
|\per A| \le \per D_{(n,r)}
\label{eq-27.16-new}\end{equation}
holds true for all values of  $n$ and $\rk A$, except the case of a nonsingular matrix of the order $4$. The equality in 
\eqref{eq-27.16-new} holds if and only if the matrix  $A$ can be obtained from  $D_{(n,r)}$ by the transposition,  row or column permutations and  multiplications of rows or columns by $(- 1)$. In the case of nonsingular   matrices  of order 4 there exists an exceptional matrix,  $D_{(4,4)}$, which is a unique exception up to the aforesaid transformations. 
\label{h1}
\end{hyp}

The confirmation of Conjecture \ref{h1}  certainly implies that Conjecture \ref{h1/1} holds. Note that  in the case $n\ge 5$ both conjectures coincide.  

Except for the bound, Conjecture \ref{h1} also claims that there exists a matrix on which the bound on the value of the permanent of $(-1,1)$-matrix of order  $n \ge 5$ and rank  $k \le n$ is achieved, i.e., it asserts that the bound is exact.

Conjecture  \ref{h1} and Problem \ref{p1} were open. They attracted considerable attention and appeared in several survey papers, see  for example  the  surveys \cite{Minc1} and \cite[Conjecture 37]{Minc-update} and references therein. 
Our paper is devoted to the complete solution of this problem. In particular we prove that the conjecture posed by Kr\"auter (Conjecture \ref{h1/1} above)  is true and moreover Conjecture \ref{h1} above is also true.

The paper is organized as follows. Section~\ref{S2} contains the main definitions and notations. In Section~\ref{S3} we consider  the case  $n\le 4$ and provide some evaluations for small $n$, which we are going to use further. Section \ref{S3.5} contains several recurrence relations for $\per D_{(n,k)}$, which are used for the induction. Section \ref{S4} contains the part of the inductive proof of the main conjecture, namely we prove that if $A\in M_{n}(\pm 1)$, $n\ge 5$, is invertible and  Conjecture~\ref{h1} holds for all $B\in M_m{(\pm1)}$, where $m<n$, then $|\per A|\le \per D_{(n,n-1)}$. In Section~\ref{S5} we investigate the behavior of the rank function on the set of submatrices of a given  $(\pm 1)$-matrix under the concatenation of  a column to a  matrix. In Section~\ref{S6} a certain vector majorization   is introduced and applied to  rank vectors of matrices under consideration. It is proved that the rank vector of $D_{(n,k,k-1)}$ is minimal among all rank vectors of matrices from $M_{k,n}(\pm 1)$. Section \ref{S7} contains the conclusion of the inductive proof, namely it is proved that if  Conjecture~\ref{h1} holds for all $B\in M_m({\pm1})$, where $m<n$, and for all $B\in GL_n({\pm1})$, then it holds for all $n\times n$ matrices. 

\section{Definitions and Notations} \label{S2}

The following notations are used through the text.
 
We say {\it signed one} for either $1$ or $-1$.
Let  $J_n$ denote the matrix of the order  $n$ consisting of ones, $O_n$ be the matrix of zeros, $I_n$ be the $n\times n$ identity matrix. We omit the index $n$ if it is clear from the context. 

For $r\ge 0$ by $\SS_{n,r}$ we denote the set of all subsets consisting of  $r$ distinct elements of the set  $\{1, \ldots, n\}$, here  $\SS_{n,0}=\emptyset$. For  $\al \in \SS_{n,r}$  by  $\ol{\al}$ we denote the set  $\{1, \ldots, n\}\setminus \al \in \SS_{n, n-r}$. 
Let $A \in M_{m,n}$, $\al \in \SS_{m, r}$, $\be \in \SS_{n,k}$. By $A(\al|\be)$ we denote the matrix obtained from the matrix  $A$ by deleting rows with numbers from $\al$ and columns with numbers from  $\be$. By $A[\al|\be]$ we denote the submatrix of $A$ located on the intersection of  the rows with numbers from $\al$ and the columns with numbers from  $\be$.

The permanent function has a well-known generalization
for the rectangular matrices $A \in M_{k,n}$ where $k \le n$:

\begin{defn}
Let  $A \in M_{k,n}$, $k\le n$. The permanent of a rectangular matrix  $A$ is defined as follows:
\begin{equation}
\per A = \suml{\al \in \SS_{n,n-k}}{} \per A(|\al).
\label{eq1}
\end{equation}
\end{defn}

In the formula \eqref{eq1} the sum is taken over the set of square submatrices of the size  $k$ belonging to the matrix $A$. So  we sum up only permanents of  square matrices. 

In the sequel we will use one more function related to the permanent   for rectangular matrices.

\begin{defn} \label{notemper}
Let $A \in M_{k,n}$ where $k \le n$. By $\mper A$ we denote the sum of absolute values of permanents for all submatrices of order $k$ of $A$. Namely,
$
\mper A = \sum_{\al \in \SS_{n,n-k}}{} |\per A(|\al)|.
$
\end{defn}
Note that in the case $A$ is a square matrix, it holds that $\mper A=|\per A|$. 

\begin{lemma}
Let  $A \in M_{k,n}$,  $k \le n$. Then 
$
|\per A| \le \mper A.
$
The equality holds if and only if all  permanents of  all square submatrices of the order  $k$ of $A$ are either simultaneously nonnegative or simultaneously nonpositive. 
\label{l14}\end{lemma}

\begin{prf}  Apply $|a+b|\le |a|+|b|$ several times.
\end{prf}

\begin{defn} \label{Stand}
Multiplication of rows or columns by $-1$, permutation of rows,  permutation of columns, transposition, and their compositions are called {\em standard\/} trans\-for\-ma\-tions.
\end{defn}

Observe that for $A \in M_n(\pm 1)$ the function  $|\per A|$ is invariant under the standard  trans\-for\-ma\-tions.

\begin{defn}
It is said that the matrices $A$ and $B$ are permutationaly equivalent if there exist permutation matrices $P,Q$ such that $A=PBQ$, i.e., $A$ can be obtained from $B$ by some permutations of rows and some permutations of columns.
\end{defn}

Thus, standard transformations 
consist of permutational equivalence, transposition and multiplication of rows or columns by $(-1)$.  
The permanent does not change under   permutation equivalence and transposition.  However, the  multiplication of a row or a column by $(-1)$ change it.

Denote the number of negative entries in a matrix $A \in M_n(\pm 1)$ by~$\nu(A)$.

\begin{lemma}
Let $A \in M_{k,n}(\pm 1)$ and  $\nu(A)=l$, where $l \le k \le n$. Assume that each row or column contains at most one negative entry. Then $A$ is permutationally equivalent to the matrix~$D_{(n,k,l)}$. 
\label{lem-1-1}\end{lemma}

\begin{prf} Permuting rows and columns of $A$ one gets $D_{(n,k,l)}$.
\end{prf}

To simplify the notations further  without loss of generality we assume that if  $A \in M_{k,n}(\pm 1)$ satisfies the conditions of Lemma~\ref{lem-1-1}, then $A$ is of the form $D_{(n,k,l)}$. Note that $\per A=\per D_{(n,k,l)}$.
 
Let $A \in M_{k,n}(\pm 1)$, $k \le n$. We denote   
$K(A) = \{A(|\al)$  for all $\al \in \SS_{n, n-k}\}.
$
 
Matrices  $C_1$  and $C_2$ from the multiset $K(A)$ are assumed to be distinct  if  they differ as the submatrices of  $A$, i.e., if they have different column indices. These matrices can coincide as matrices in  $M_k(\pm 1)$. For example, all matrices in  $K(J_{k,n})$ are evidently equal to $J_k$, however we  consider them as separate elements. Note that for any $A \in M_{k,n}$ we have $|K(A)|={n\choose{k}}$.

The permanent of   $A \in M_{k,n}$ can be rewritten as
$
\per A = \sum_{B\in K(A)}{} \per B.
$

Below we suggest the rank stratification of the set $K(A)$. 

\begin{defn} \label{def_rk_vect}
Consider the multiset  $X$ consisting of some, possibly repeating, matrices from $M_k(\pm 1)$. The {\em rank vector\/}  $R(X) = (r_1,\ldots, r_k) \in \N^k$ of   $X$ is defined by 
$
r_i = |\{A \in X| \rk A = k-i+1\}|.
$ 
\end{defn}

Observe that the first component of the vector  $R(X)$ 
is equal to the number of matrices of  rank   $k$ (i.e., invertible matrices)  in the multiset 
 $X$, the second component of  $R(X)$ is equal to the number of matrices of rank  $k-1$ in $X$, and so on.
The last component  $r_k$  is equal to the number of rank one matrices in  $X$. Since  $X$ consists of $(-1,1)$-matrices, we can conclude that the minimal rank of any matrix from $X$ is equal to $1$. It follows that the sum of all components of the vector  $R(X)$ is equal to the number of matrices in~$X$.

\begin{defn}
Let $A \in M_{k,n}(\pm 1)$ where $k \le n$. The rank vector of the matrix  $A$ is defined to be  $R(A) = R(K(A))$.
\end{defn}

\begin{defn} \label{defCb}
Let $C \in M_k(\pm 1)$ and $b$ be a vector of size $k$ with elements from $\{-1,1\}$. Denote by  $<C, b>$ the multiset of matrices  $\{C_1, \ldots, C_k\}$ where each  $C_i$ is obtained by concatenation of the matrix $C$ and the column $b$ and  deleting the $i$th column of    $C$. Here we include repeatedly the equal matrices, if necessary. Below we say that column $b$ replaces $i$th column in the matrix~$C$.  
\label{t-08-1}\end{defn}

Further we need the following technical lemma.

\begin{lemma}
Let $A \in M_n(\pm 1)$  be a nonsingular matrix. 
Then by the standard transformations   $A$ 
can be reduced to a matrix  $B$, satisfying the following conditions:
\begin{enumerate}
	\item The first row and the first column of   $B$ consist of the entries $1$ only.
	
	\item The submatrix $B(1|1)$ has rank  $(n-1)$.
	
	\item  $|\per A| = |\per B|$.
\end{enumerate}
\label{lem1}
\end{lemma}

\begin{prf}
Since  $A$ is a nonsingular matrix, there exists nonsingular its submatrix
of order  $(n-1)$. By the permutation of rows and columns of matrix $A$ we obtain $A'=(a_{ij}')$ such that its submatrix  $A'(1|1)$ is nonsingular. 

For $i = 1, \ldots, n$ such that $a_{1i}' = -1$ we multiply the $i$th column by $-1$ to obtain $A''=(a_{ij}'')$. For $j = 2, \ldots, n$ such that $a_{j1}'' = -1$ we multiply the $j$th row by $-1$ to obtain $A'''=B$. Multiplication of rows and columns by $-1$ does not change the rank. Thus $B(1|1)$ is nonsinguar, which proves $1.$ and $2$.

Since standard transformations preserve the absolute value  of the permanent, 3. is also true.
\end{prf}

\section{Matrices of small sizes}  \label{S3}

Our proof  is inductive. 
So let us collect in this section the  results for matrices of small size. We are going to prove  that the conjecture is true for  $n = 2$ and $3$  and to investigate the situation for $n=4$.
Note that for $n = 4$ there exists an exceptional example showing that Conjecture  \ref{h1/1} is not true. 

The following lemma is a straightforward consequence of Lemma  \ref{lem1}. 
\begin{lemma} \label{lem_n2}
Let $A\in M_2(\pm 1)$ be invertible. Then $|\per A| = \per D_{(2,1)} = 0$.
\end{lemma}

\begin{lemma} \label{lem_n3}
Let $A\in M_3(\pm 1)$ and $\rk A = k$. Then $|\per A| = \per D_{(3,k-1)}$ and $A$ can be reduced to $D_{(3,k-1)}$ by standard transformations.
\end{lemma}

\begin{prf} 
While in $A$ there is a row or a column with at least two $-1$s let us multiply this row or column by $-1$. After each   multiplication the number of negative entries in $A$ strictly decreases. Thus we obtain a matrix $A'$ with no two $-1$s in any row or column. By Lemma \ref{lem-1-1} $A'$ is equivalent to $A'' = D_{(3,i)}$ for some $i\le 3$.

If $i < 3$, then $\rk A' = i+1$ and the lemma is proved. If $i = 3$, then we multiply the first row, and after that the second and third columns by $-1$. There are two $-1$s in the obtained matrix $A''$. They are  $a_{23}'' = a_{32}'' = -1$. By Lemma \ref{lem-1-1} $A''$ is equivalent to $D_{(3,2)}$. The lemma is proved.
\end{prf}

\begin{cons} \label{Cor3.3}
  Conjecture   $\ref{h1}$ is true  for $n=2,3$. 
\end{cons}

\begin{prf}
Lemma  \ref{lem_n2} proves Conjecture   $\ref{h1}$  for invertible matrices with  $n=2$.
It is straightforward to see that all singular $(-1,1)$ matrices of size 2 can be reduced to $J_2$ by standard transformations.
Lemma  \ref{lem_n3} proves Conjecture   $\ref{h1}$  in the case  $n=3$. 
\end{prf}

We need also the following examples of concrete evaluation of the permanent function for several matrices of special type.

\begin{lemma} \label{lem:n4}
The case  $n = 4$ is exceptional. Conjecture \ref{h1/1} does not hold in this case.
\end{lemma}

\begin{prf}
Expanding  
on the last row we have  $
\per D_{(4,3)} =
 3\per D_{(3,2)} + \per D_{(3,3)} = 6 - 2 = 4.
$
 Then $\per D_{(4,4)} = 8 > 4 = \per D_{(4,3)}$.
\end{prf}

Below we recall the generalized Laplace formula on the permanent expansion on  several rows  
 which will be useful in the sequel.

\begin{theorem}[Generalized Laplace formula, {\cite[Chapter 2.1, Theorem 1.2]{Minc}}]
Let $A \in M_n$ and  $\be \in \SS_{n,k}$, $0 < k < n$. The following equality is true:
$$
\per A = \suml{\al \in \SS_{n,k}}{} \per A(\al|\be) \per A[\al|\be].
$$
\end{theorem}

\begin{lemma}
Let $A \in M_4(\pm 1)$. Then  $\per A$ is divisible by $4$.
\label{del}\end{lemma}

\begin{prf}
Observe that the permanent of $(\pm1)$-matrix of order $2$ is divisible by $2$. Let us apply Laplace expansion formula by the first two rows  to the matrix of order $4$. Then we obtain the sum  with the property that 
 any summand is a product of even numbers. Therefore, it is divided by $4$. 
 This implies the statement.
\end{prf}

Now we classify all invertible $(\pm1)$-matrices of the size $4$, for which Conjecture~\ref{h1/1} does not hold. As it will be shown further, Conjecture~\ref{h1/1} holds for singular matrices of the size $4$, so we will describe the unique possible case when Conjecture~\ref{h1/1} cannot be generalized. 

\begin{theorem} \label{lem_n4}
Let a matrix $A \in GL_4(\pm 1)$ be such that  $|\per A| = \max_{C \in GL_4(\pm 1)}|\per C| $. Then $|\per A| =8$ and using the standard transformations we can reduce  $A$ to the matrix $D_{(4,4)}$. If for $B  \in GL_4(\pm 1)$ it holds that $|\per B| < 8$, then  $|\per B| \le \per D_{(4,3)} = 4$. 
\label{experiment}
\end{theorem}

\begin{prf}
Let  $A \in M_4(\pm 1)$ be such that $|\per A|$ is maximal. By Lemma \ref{lem1} (item $1$) the matrix $A$ can be reduced to the form with $a_{1i} = 1$, $i = 1, \ldots, 4$.

If there is a row in $A$ with at least three $-1$s, then we multiply this row by $-1$. Thus we may consider only the case when there are at most two $-1$s in each row of~$A$. 

{\bf 1.} Suppose  there is at most one $-1$ in each row of $A$. Then $A$ has at most three $-1$s. If there are two negative entries in the same row (or column), then there are two columns (or rows) with all positive entries. Hence, $A$ is singular. Thus by Lemma \ref{lem-1-1} $A$ is permutationally equivalent to $D_{(4,3)}$.    

{\bf 2.} Suppose there is a row with two $-1$ in $A$ . Up to the permutation of rows and columns we may assume
$ A =
\left(
\begin{smallmatrix}
1 & 1 & 1 & 1\\
1 & 1 & -1 & -1 \\
a_{31} & a_{32} & a_{33} & a_{34} \\
a_{41} & a_{42} & a_{43} & a_{44}
\end{smallmatrix} \right).
$

For each $i \in \{3,4\}$ such that $a_{i1} = -1$ we multiply $i$th row by $-1$. By the generalized Laplace formula applied to first $2$ rows we get
\begin{equation}
|\per A| = |2\per
\left( \begin{smallmatrix}
a_{33} & a_{34} \\
a_{43} & a_{44}
\end{smallmatrix} \right)
- 
2\per
\left( \begin{smallmatrix}
1 & a_{32} \\
1 & a_{42}
\end{smallmatrix} \right)| \le 4+4=8\label{eq:n10}
\end{equation}

Since there is no proportional columns in $A$ and $|\per A|$ is maximal, we obtain two consequences: 

{\bf A.} $|\per
\left( \begin{smallmatrix}
1 & a_{32} \\
1 & a_{42}
\end{smallmatrix} \right)| = 2$ and $a_{32} = a_{42} = -1$. 

{\bf B.} $|\per
\left( \begin{smallmatrix}
a_{33} & a_{34} \\
a_{43} & a_{44}
\end{smallmatrix} \right)| = 2$ and there is at least one $1$ in submatrix $A(1,2|1,2)$. Since $\rk A(1,2|1,2) = 1$ and columns are proportional and columns of $A$ are not proportional, we obtain $a_{44} = 1$ and $a_{34} = a_{43} = -1$.   

Hence,
$
A = 
\left(\begin{smallmatrix}
1 & 1 & 1 & 1\\
1 & 1 & -1 & -1 \\
1 & -1 & 1 & -1 \\
1 & -1 & -1 & 1
\end{smallmatrix} \right).
$

Now, multiplying the rows $2$, $3$, $4$ and the column $1$ by  $(-1)$ we obtain the matrix $D_{(4,4)}$. So, if the maximum of absolute value of permanent is achieved on $A$, then $A$ can be reduced to $D_{(4,4)}$ by standard transformations. 
\end{prf}

We conclude by several explicit computations which we need later on.

\begin{proposition}\label{PropSmallD}
The following equalities hold: \\
$
\begin{smallmatrix}
\begin{cases}
\per D_{(4,0)} = 24 \\
\per D_{(4,1)} = 12 \\
\per D_{(4,2)} = 8 \\
\per D_{(4,3)} = 4 \\
\per D_{(4,4)} = 8
\end{cases} & \ \ \ & \ \ \ & \ \ \ & \ \ \ & \begin{cases}
\per D_{(6,6)} = 112\\
\per D_{(5,5)} = 8\\
\per D_{(4,4)} = 8\\
\per D_{(3,3)} = -2\\
\per D_{(2,2)} = 2\\
\per D_{(1,1)} = -1
\end{cases} \end{smallmatrix}\enspace .
$
 \end{proposition}

Since  $\per D_{(3,3)} = -2$, we have that in the expansion for permanent of the  matrix $D_{(4,3)}$ by the last row there is a negative summand. Hence $\per D_{(4,3)}=4<8=\per D_{(4,4)}$. In the case  $n > 4$, there are no negative summands in the expansion  of  $\per D_{(n,n-1)}$ by the last row.

In particular this explains the main reason why in  the case  $n = 4$ the conjecture does not hold.
More details about this  case   can be found in~{\cite[Lemma 2 and Remark 2]{KrSe2}}. 

\section{Several recurrent relations for $\per D_{(n,k)}$}  \label{S3.5}

Let us remind the following lemma which is proved in \cite{KrSe}.

\begin{lemma}  \label{lem40}
Let $n \ge 5$. Then the following statements are true:
\begin{enumerate}
		\item $\per D_{(n,n)} > 0$, see \cite[Lemma $2$, Formula $(4)$]{KrSe}. 
	\item $\per D_{(n,k)} > 0$ for all~$k$, see  \cite[Lemma $2$, Formula $(4)$]{KrSe}. 
\item $\per D_{(n,k-1)} > \per D_{(n,k)}$ for all~$k$, see \cite[Lemma $2$, Formula $(5)$]{KrSe}.  
	\end{enumerate}
 \end{lemma}

Lemma  \ref{lem40} has the following corollary:

\begin{cons}
Let $n \ge 5$. Then for any  $l, k$, $0 \le l < k \le n$, we have $\per D_{(n, l)} > \per D_{(n, k)}$.
\label{con40}\end{cons}

We prove here the identity which can be also derived from the identities  \cite[Lemma 1, Formulas (2) and (3)]{KrSe}. However, for the later use we need the identity in the form presented below.
\begin{lemma} \label{lem:c}
$
\per D_{(n,k-1)} = \per D_{(n,k)} + 2\per D_{(n-1,k-1)}
$ for all $k\le n$. 
\end{lemma}

\begin{prf}
Let us decompose  permanent of the matrix  $D_{(n,k-1)}$ by any row which does not contain $-1$s. 
This decomposition  is possible  since  $n > k-1$:

$$
\per D_{(n,k-1)} = \suml{i = 1}{k-1} \per D_{(n-1,k-2)} + \suml{i = k}{n} \per D_{(n-1, k-1)} = 
$$
$$
(k-1)\per D_{(n-1, k-2)} + (n-k+1)\per D_{(n-1,k-1)}.
$$
Now decompose   $\per D_{(n,k)}$ by the first row which contains one $-1$:
$$
\per D_{(n,k)} =
 -\per D_{(n-1,k-1)} + \suml{i = 2}{k} \per D_{(n-1, k-2)} + \suml{i = k+1}{n} \per D_{(n-1,k-1)} = 
$$
$$
(k-1)\per D_{(n-1,k-2)} + (n-k-1)\per D_{(n-1,k-1)}
$$
By subtracting we get
$
\per D_{(n,k-1)} - \per D_{(n,k)} = 2\per D_{(n-1, k-1)}. 
$
\end{prf}

We are going to obtain a recurrent formula for $\per D_{(n,n)}$ via $\per D_{(k,k)}$. Note that the following lemma can also be   derived from  \cite[Lemma 1]{KrSe}.

\begin{lemma}
Let $n \ge 5$. The following identity holds: 
$$
\per D_{(n,n)} = (n-2)\per D_{(n-1,n-1)} + (2n-2)\per D_{(n-2,n-2)}.
$$
\label{lem:c2}
\end{lemma}

\begin{prf}
Let us expand the permanent of  $D_{(n,n)}$ by the first row:
\begin{equation}
\per D_{(n,n)} = -\per D_{(n-1, n-1)} + (n-1)\per D_{(n-1, n-2)}.
\label{eq5}
\end{equation}
By Lemma \ref{lem:c} we have the equality: 
\begin{equation}
\per D_{(n-1,n-1)} = \per D_{(n-1,n-2)} - 2\,\per D_{(n-2,n-2)}.
\label{eq6}
\end{equation}
We substitute the expression \eqref{eq6} into the expression  \eqref{eq5} and get
\begin{equation}
\per D_{(n,n)} = (n-2)\per D_{(n-1,n-2)} + 2\,\per D_{(n-2,n-2)}.
\label{eq7}
\end{equation}
Applying Lemma \ref{lem:c} to the permanent of   $D_{(n-1,n-2)}$ we obtain the required identity:
$
\per D_{(n,n)} = (n-2)\per D_{(n-1,n-1)} + (2n-2)\per D_{(n-2,n-2)}.
$
\end{prf}

\section{The proof of the Kr\"auter conjecture for nonsingular matrices}  \label{S4}

In this section we prove Conjecture \ref{h1} for nonsingular matrices of order $n \ge 5$. This is a part of inductive proof and we suppose that Conjecture \ref{h1} holds for all matrices of order less than $n$. 

\subsection{Permanent preserving reductions of nonsingular matrices}

In this subsection we construct the reduction for nonsingular matrices of size $n \ge 5$ to the specific form which is useful for estimate of permanent of a matrix. The reduction is made by standard transformations thus it preserves absolute values of permanent. We also show that some matrices of order $5$ and $6$ cannot be reduced to required form. We characterize all such matrices up to standard transformations.  

\begin{denote}
By $Q_m \in M_m(\pm 1)$ we  denote the following matrix $Q=(q_{ij})$:
$$
q_{ij} = \begin{cases}
-1, \mbox{  if } i = j;\\
-1,  \mbox{  if } j = i+1; \\
-1,  \mbox{  if }  i = m  \mbox{ and  } j = 1;\\
1, \mbox{ otherwise.}
\end{cases}
$$
In other words, $Q_m=J_m-2I_m-2(E_{12}+E_{23}+\ldots+E_{n-1,n}+E_{n,1})$.
\end{denote}

\begin{exam}
 $
Q_2 =  \left(
\begin{smallmatrix}
-1 & -1 \\
-1 & -1
\end{smallmatrix} \right),  $
$Q_3 = \left(
\begin{smallmatrix}
-1 & -1 & 1 \\
1 & -1 & -1 \\
-1 & 1 & -1
\end{smallmatrix} \right),  $ $
Q_4 = \left(
\begin{smallmatrix}
-1 & -1 & 1 &1\\
1 & -1 & -1 & 1\\
1 & 1 & -1 & -1\\
-1 & 1 &  1 & -1
\end{smallmatrix}\right).  
$
\end{exam}
It is straightforward to see that for any positive $k$ the matrix $Q_{2k}$ is singular and   the matrix $Q_{2k + 1}$ is nonsingular. 

\begin{lemma} 
Let $A \in M_n(\pm 1)$ and each row and column of $A$ contain exactly two negative entries. Then by a permutations of rows and columns $A$ can be reduced to the form 
\begin{enumerate}
	\item there are blocks  $Q_{k_i}$ on the diagonal of $A$, $k_i\ge 2$, $\sum k_i = n$;
	
	\item all other entries are equal to $1$. 
\end{enumerate}
\label{lem-1-3}\end{lemma}

\begin{prf} Let us choose the appropriate permutations subsequently.

{\bf 1.} Since there are two $-1$s in the first row of $A$, up to a permutation of columns, we may assume that $a_{11} = a_{12} = -1$.

{\bf 2.} Since there are two $-1$s in the second column of $A$ and one of them lies in the first row, we can without loss of generality assume  $a_{12} = a_{22} = -1$, permuting, if necessary, the second row with the row containing the other~$-1$. 

{\bf 3.} Then for any row $k\ge 2$ we do the following permutations of rows and columns. By the previous steps we have  $a_{k-1,k-1} = a_{k-1,k} = a_{k,k} = -1$. 

{\bf 3.1.} If there exists $l>k$ with  $a_{k,l} = -1$, then we permute the columns with the indices $l$ and $(k+1)$. Due to the location of negative entries in the rows with the indices $1,\ldots, k-1$, we have  $a_{m,k+1} = 1$ for all $m < k$. Thus, there exists $p>k$ for which  $a_{p,k+1} = -1$. Then after permuting  the $p$th and  $(k+1)$st rows we get  $a_{k,k} = a_{k, k+1} = a_{k+1,k+1} = -1$. So, we can apply again Step {\bf 3} for the next value of $k$.

{\bf 3.2.} If there is no $l>k$ such that  $a_{k,l} = -1$, then $a_{k,1}=-1$, since each row and column contains exactly two negative entries.

{\bf 4.} We repeat Step {\bf 3} until {\bf 3.2} does not hold. If {\bf 3.2} holds, then  $A[1,\ldots, k| 1,\ldots, k] = Q_{k}$. Note that all entries of $A[1,\ldots, k|k+ 1,\ldots, n]$ and $A[k+1,\ldots, n| 1,\ldots, k]$ are equal to 1 since only two $-1$s can be in a row. This guarantees the condition 2. of the lemma. Since each row and column of $Q_k$ contains two $-1$s, the conditions of the lemma hold for the matrix  $A'=A(1,\ldots, k| 1,\ldots, k)$ as well. Repeating Steps {\bf 1} --- {\bf 3} for $A'$ we find the next block  $Q_{k'}$. Then repeating this process and noting that  the size of $A'$ decreases monotonically, we obtain a required form. 
\end{prf}

Let $P_1, P_2 \in M_6(\pm 1)$ denote the following matrices: \\
$
P_1 = \left(
\begin{smallmatrix}
	1 & 1 & 1 & 1 & 1 & 1 \\
	1 & -1 & -1 & 1 & 1 & 1 \\
	1 & 1 & -1 & -1 & 1 & 1 \\
	1 & 1 & 1 & -1 & -1 & 1 \\
	1 & 1 & 1 & 1 & -1 & -1 \\
	1 & -1 & 1 & 1 & 1 & -1 \\
\end{smallmatrix}\right)$, 
$  
P_2 = 
\left(\begin{smallmatrix}
1 & 1 & 1 & 1 & 1 & 1\\
-1 & -1 & 1 & 1 & 1 & 1\\
-1 & 1 & -1 & 1 & 1 & 1\\
-1 & 1 & 1 & -1 & 1 & 1\\
1 & 1 & 1 & -1 & -1 & 1\\
1 & 1 & 1 & -1 & 1 & -1\\
\end{smallmatrix} \right).
$ \\
Direct computations show that $\per P_1 =  \per P_2 = 16$.

Now we formulate and prove the following corollary from Lemma~\ref{lem-1-3}.

\begin{cons}
Let $A \in M_6(\pm 1)$ be an arbitrary nonsingular matrix satisfying  $a_{1,k} = a_{k,1} = 1$ for all $k = 1, \ldots, 6$. Assume also that each row and each column of $A$ from second till $6$th contains exactly $2$ negative entries. Then $A$ is permutationally equivalent to~$P_1$.
\label{con-1-1}
\end{cons}

\begin{prf}
We apply Lemma \ref{lem-1-3} to the submatrix $A(1|1)$, permuting the rows and the columns of $A$. Since the first row and column of $A$ consist only of $1$s, after permutations they are unchanged. 

By Lemma \ref{lem-1-3} and since the size of each block $Q_j$ is greater than or equal to $2$, the submatrix $A(1|1)$ can be transformed to one of the following forms.

{\bf 1.} There are blocks  $Q_2$ and $Q_3$ on the diagonal, and the rest of the entries are equal to $1$. In this case the rows of $A$ corresponding to the block $Q_2$ are equal. This contradicts the invertibility of $A$. So, this case does not hold.

{\bf 2.} There is a block $Q_5$ on the diagonal. Hence, $A$ is $P_1$.
\end{prf}

\begin{defn}
Let $A$ be a $(-1,1)$-matrix. We say that  $A$ satisfies the {\em condition ${\mathfrak A}$\/}  if  the following two facts are true:
\begin{enumerate}
	\item all entries of the first row of $A$ are $1$,	
	\item the second row of $A$ contains at least three positive entries and at least three negative entries. 
\end{enumerate}
\end{defn}

The next lemma  allows to reduce a nonsingular matrix by standard transformations to the form in
which we can efficiently estimate the permanent function. 

\begin{lemma}
Let $A\in M_n(\pm 1)$ be nonsingular, $n\ge 6$. 
Then $A$ can be reduced by the standard transformations to one of the following forms: \\
1. $D_{(n,n-1)}$. \\
2. $D_{(n,n)}$. \\
3. $P_1$ or $P_2$. \\
4. A matrix satisfying the condition ${\mathfrak A}$.
\label{lem:vid}
\end{lemma}

\begin{prf}
{\bf 1.} By Lemma \ref{lem1} we can reduce $A$ by the standard transformation to the matrix $A'$ such that all entries in the first row are $1$. Without loss of generality we shall further denote $A'$ as $A$. Suppose, there is an index $i$ such that $i$th row of $A$ contains at least three $1$s and at least three $-1$s. We can permute the $i$th row and the second row of the matrix. Obtained matrix satisfies condition ${\mathfrak A}$ and the lemma is proved in this case. 

{\bf 2.} So, further we suppose that any row of $A$ contains either at most two positive or at most two negative entries. Equivalently, each row contains either at least $n-2$ entries $1$ or   at least $n-2$ entries~$-1$.

{\bf 3.} We multiply by $-1$ all  rows of $A$ which contain more than two $-1$s. The obtained matrix we denote again  by $A$ to simplify the notations. Therefore we can assume in the sequel that $A$ satisfies the following conditions:

{\bf C1.} All entries in the first row of $A$ are $1$. 

{\bf C2.} Each row of $A$ contains no more than two entries $-1$.

{\bf C3.} $A$ is invertible.

It follows from the invertibility of $A$,  that only the entries in the first row are all positive. Below we consider different cases of the location of negative entries in the rows of $A$.   

{\bf 4.} Assume firstly that each row  of $A$, except the first, contains exactly one $-1$. Then $\nu(A) = n-1$. If there is a column of $A$ with at least two $-1$s, then by the pigeonhole principle there are at least two columns without $-1$s. These   contradicts  the invertibility of $A$. Then  $A$ can be reduced to $D_{(n, n-1)}$ by Lemma \ref{lem-1-1}.  

{\bf 5.} Let us denote  by $k$ the number of the rows of $A$ with one $-1$. Hence, $k<n$ by Condition {\bf C1}. The case $k=n-1$ is considered in Item {\bf 4}. So, we now consider the case $k<n-1$. By Conditions {\bf C1} and {\bf C2}, there are $(n-k-1)$ rows of $A$ with two $-1$s and $\nu(A) = k+2(n-k-1)=2n-k-2$. 
We can permute the rows in $A$ in such a way that each of the  rows with the indices $2, \ldots, k+1$ contains one $-1$, and each of the rows with the indices $k+2, \ldots, n$ contains two  $-1$s. 

{\bf 6.} For each $i$ such that $2 \le i \le k+1$ we can permute the columns of $A$ in order to obtain $a_{i, i} = -1$. Indeed, $A$ is  nonsingular, so no two $-1$s in the rows $2, \ldots, k+1$ lie in the same column. Without loss of generality we denote the obtained matrix again by $A$. Thus from now on we have the equality $A[2, \ldots, k+1|2, \ldots, k+1] = J_k - 2I_k$.
We proceed further by separate considerations of the cases $k\ge 3, \ k=2, \ k=1, \ k=0$. 

{\bf 7.} Let $k \ge 3$. There are only two entries $-1$ in the row with the index $(k+2)$. Hence by the pigeonhole principle there is $i \in \{2,\ldots, k+1\}$ such that $a_{k+2,i} = 1$. Note that by Item {\bf 6} the entry $a_{i,i}=-1$. Let us multiply the $i$th column by $-1$. In obtained matrix $A'$ there are three $-1$s in $(k+2)$nd row and all entries in the $i$th row are $1$s. By permutations of the first row with the $i$th row and the second row with the $(k+2)$nd, we obtain a matrix satisfying the condition ${\mathfrak A}$. Thus the lemma is proved in this case.

{\bf 8.} Let $k = 2$. We can find   $i > 3$ and $j = 2$ or $3$ such that $a_{ij} = 1$. Indeed, if there is no such $ {i,j}$, then both $-1$s of each of rows with the indices $4, \ldots, n$ lie in the second and in the third columns. Hence all rows with the indices $4, \ldots, n$ are equal to each other and $A$ is singular since $n\ge 6$. This contradicts   {\bf C3}.

Thus multiplying $j$th column by $-1$ and  permuting the $j$th row with the first row, and the $i$th row with the second row we get the matrix satisfying condition ${\mathfrak A}$, so the lemma is proved in this case also. 

{\bf 9.} Let $k = 1$. 
If there is $i>2$ such that  $a_{i2} = 1$, then by multiplication of the second column by $-1$ and application of the cycle permutation $(2,1,i)$ to the second, first, and $i$th  rows, we obtain a matrix satisfying   ${\mathfrak A}$. 

Suppose that  $a_{i2}=-1$ for all $i>2$. It means that there are $(n-1)$ negative entries in the second column (since $a_{22}=-1$ by Item {\bf 6}) and $\nu(A) = 2n-3$ by Item {\bf 5}. 
Let us multiply the second column by $-1$ to obtain $A'$.  Then $\nu(A') = n-1$. By repeating the arguments from Item {\bf 4} we obtain that  $A'$ can be reduced to  $D_{(n,n-1)}$ by standard transformations, thus $A$ can be reduced. 

{\bf 10.} Let $k = 0$. It means that each row of $A$ except the first row contains exactly two $-1$s. Hence $\nu(A) = 2n-2$. 
Let us consider several cases now. These cases depend on the quantity of $-1$s in columns of $A$. Note that each column contains at most  $n-1$ negative entries since the first row is positive. 

{\bf 10.1.} Suppose that there is a column in $A$ with $(n-1)$ negative entries, i.e., $a_{ij} = -1$ for some $j$ and all $i = 2, \ldots, n$. Multiplying the $j$th column of $A$ by $-1$ we obtain $A'$ such that $\nu(A') = 2n-2-(n-1)+1=n$. 
Also each row of $A'$ contains exactly one $-1$. Since $A'$ is nonsingular, each column 
also contains exactly one $-1$. Thus $A'$ is permutationally equivalent to $D_{(n,n)}$, so $A$ can be reduced to $D_{(n,n)}$ by standard transformations.

{\bf 10.2.} Suppose that there is a column in $A$ with $m \in \{3, \ldots, n-2\}$ negative entries. Without loss of generality $a_{i1} = -1$ where $i = 2, \ldots, m+1$. Since $A$ is nonsingular and there are two $-1$ in each row $2, \ldots, n$, we can permute columns $2, \ldots, n$ to obtain $a_{ii} = -1$ for $i = 2, \ldots, m+1$. 

We fix $i,j$ such that $a_{m+2,i} = a_{m+2,j} = -1$. Due to the structure of $A$ since $m \ge 3$, there is an index $l \in \{2, \ldots, m+1\}\setminus \{i,j\}$  such that $a_{l,1} = a_{l,l} = -1$. There exist two subcases depending on $n$. 

{\bf 10.2.1.} Suppose that $n \ge 7$. We multiply the columns $i,j$ by $-1$. Then all entries of the row  $(m+2)$ are positive and there are four $-1$s  in the $l$th row, namely,   $a_{l,1}= a_{l,l}= a_{l, i}= a_{l,j}=-1$. Hence there are  $n-4$ positive entries in the $l$th row. Since $n\ge 7$, we obtain $n-4\ge 3$. Therefore, permuting the first row with the $(m+2)$nd row and the $l$th row with the second row we obtain a matrix satisfying condition ${\mathfrak A}$ and the lemma is proved in this case.

{\bf 10.2.2.} It remains to consider the case $n = 6$. Recall that by the conditions of the Item {\bf 10.2} we have $m \in \{3, \ldots, n-2\}=\{3,4\}$.


{\bf Case $m = 3$}. Then there are three $1$s in  the first column. If there is also a column of all $1$s, then $A^t$ contains a row with all elements equal to $1$ and a row with three $1$s and three $-1$s. We can permute rows of $A^t$ in the way that the obtained matrix satisfies the condition ${\mathfrak A}$.

So, assume now that there is no columns in $A$ with all $1$s. By Item {\bf 10.2} 
$
A(1,5,6|) = 
\left(\begin{smallmatrix}
	-1 & -1 & 1 & 1 & 1 & 1\\
	-1 & 1 & -1 & 1 & 1 & 1\\
	-1 & 1 & 1 & -1 & 1 & 1
\end{smallmatrix} \right).
$

Since there are no columns with all $1$s, both columns with indices $5$ and $6$ contain the entries $-1$. These entries cannot be in the same rows since there are just $2$ negative entries in each row by the conditions of Item {\bf 10} and $A$ is invertible by assumptions. Hence up to permutations of rows and columns  only the following $2$ options are possible for   $A[5,6|5,6]$:
$ 
\left(\begin{smallmatrix}
-1 & 1 \\
1 & -1
\end{smallmatrix} \right)
 \text{ or }
\left(\begin{smallmatrix}
-1 & 1 \\
-1 & -1
\end{smallmatrix} \right).
$

In the last case the other $-1$s in the $5$th row can be in the $i$th column, where $i\in\{2,3,4\}$. Permuting the $i$th and $4$th columns and then the $i$th and $4$th rows we can consider $A$ in the form 
$
A_1 = 
\left(\begin{smallmatrix}
1 & 1 & 1 & 1 & 1 & 1\\
-1 & -1 & 1 & 1 & 1 & 1\\
-1 & 1 & -1 & 1 & 1 & 1\\
-1 & 1 & 1 & -1 & 1 & 1\\
1 & 1 & 1 & -1 & -1 & 1\\
1 & 1 & 1 & 1 & -1 & -1\\
\end{smallmatrix} \right).
$

In the first case the remaining two negative entries can be either in one column or in different columns. Up to the permutations as above the following  possibilities occur:
$
A_2 = 
\left(\begin{smallmatrix}
1 & 1 & 1 & 1 & 1 & 1\\
-1 & -1 & 1 & 1 & 1 & 1\\
-1 & 1 & -1 & 1 & 1 & 1\\
-1 & 1 & 1 & -1 & 1 & 1\\
1 & 1 & 1 & -1 & -1 & 1\\
1 & 1 & 1 & -1 & 1 & -1\\
\end{smallmatrix} \right), \ 
A_3 = 
\left(\begin{smallmatrix}
1 & 1 & 1 & 1 & 1 & 1\\
-1 & -1 & 1 & 1 & 1 & 1\\
-1 & 1 & -1 & 1 & 1 & 1\\
-1 & 1 & 1 & -1 & 1 & 1\\
1 & 1 & -1 & 1 & -1 & 1\\
1 & 1 & 1 & -1 & 1 & -1\\
\end{smallmatrix} \right).
$

Note that $A_2 = P_2$. We multiply the last row of  matrices $A_1$ and $A_3$ by $-1$ and transpose them. Both obtained matrices    contain a row with all $1$s and a row with three $1$s and three $-1$s.  Thus the condition ${\mathfrak A}$ is satisfied.

{\bf Case $m=4$.} Depending on whether $a_{66} = 1$ or $a_{66} = -1$   $A$ is equal to 
$
A_4 = 
\left(\begin{smallmatrix}
1 & 1 & 1 & 1 & 1 & 1\\
-1 & -1 & 1 & 1 & 1 & 1\\
-1 & 1 & -1 & 1 & 1 & 1\\
-1 & 1 & 1 & -1 & 1 & 1\\
-1 & 1 & 1 & 1 & -1 & 1\\
1 & 1 & 1 & 1 & -1 & -1\\
\end{smallmatrix} \right) $ or 
$A_5 = 
\left(\begin{smallmatrix}
1 & 1 & 1 & 1 & 1 & 1\\
-1 & -1 & 1 & 1 & 1 & 1\\
-1 & 1 & -1 & 1 & 1 & 1\\
-1 & 1 & 1 & -1 & 1 & 1\\
-1 & 1 & 1 & 1 & -1 & 1\\
1 & 1 & 1 & -1 & -1 & 1\\
\end{smallmatrix} \right)
$ up to permutations of rows and columns.   
For matrices $A_4$ and $A_5$ we multiply the second row by $-1$ and transpose the matrices. Then   the condition ${\mathfrak A}$ is satisfied.

{\bf 10.3.} Suppose now that there are at most two $-1$s in each column of $A$ and there is a column with all $1$s. Say, all elements of the first column are positive. Then by Lemma \ref{lem-1-3} we can permute rows and columns such that the submatrix $A(1|1)$ contains several diagonal blocks $Q_{m_i}$ and all other entries of $A$ are equal to $1$. 

Let $n \ge 7$. If there is a block $Q_{n-1}$ in $A$, we choose the indices $j$ and $k$ corresponding to the first and third its rows. If there is no block $Q_{n-1}$ in $A$, then there are at least two  blocks $Q_{m_i}$ in $A$, and we choose $j$ and $k$ in such a way that $j$th and $k$th rows of $A$ lie in the different blocks. In both cases all  four $-1$s  lying in the rows $j$ and $k$ are in different columns. Say, they are $a_{js}, a_{jt}, a_{ku},a_{kv}$, where $s,t,u,v$ are all different. Then after multiplication of $u$th and $v$th columns with $-1$ we obtain that $j$th row of the obtained matrix contains four $-1$s (and hence, at least three 1s since $n\ge 7$), and $k$th row is positive. Hence after the permutation of the $k$th and first rows and the $j$th and second rows we obtain a matrix satisfying the  condition ${\mathfrak A}$. 

In the case $n = 6$, by Corollary \ref{con-1-1} we obtain the matrix $P_1$ and the lemma is proved.

{\bf 10.4.} Suppose that each column of $A$ contains either one or two $-1$s. Since  $\nu(A) = 2n-2$, there are two columns with one $-1$ in $A$. Let $a_{21} = a_{22} = -1$ and $a_{21}$ is the unique negative entry in its column. We multiply the second row of $A$ by $-1$ to obtain $A'$. All entries of the first column of $A'$ are $1$. Since $n \ge 6$, there is a column, say $i$th, in $A$ containing two $1$s and satisfying $a_{2i} = 1$. It means that in $A'$ in the $i$th column there are $3$ negative entries and $n-3 \ge 3$ positive entries. Let us permute the second and the $i$th columns of $A'$ and transpose the matrix. Obtained matrix satisfies the condition ${\mathfrak A}$. Thus the lemma is  finally proved.    
\end{prf}

In the following lemma we prove an analog of Lemma \ref{lem:vid} for  $n = 5$. 

\begin{lemma}
Let $A \in M_5(\pm 1)$ be  nonsingular. Then we can reduce  $A$ by standard transformations to one of the following forms:
 \\
1. $D_{(5,4)}$, \\
2. $B=(b_{ij})\in M_5(\pm 1)$, where $b_{1,k}=1$ for $k = 1, \ldots, 5$, $b_{2,1} = b_{2,2} = -1$ and $b_{2,3} = b_{2,4} = b_{2,5} = 1$.
\label{d5}
\end{lemma} 

\begin{prf}

{\bf 1}. Multiply the columns of $A$ by $1$ or $-1$ to obtain $A'=(a'_{ij})$ satisfying $a'_{1,k} = 1$ for $k = 1, \ldots, 5$. To simplify the notations we  further assume $A=A'$.

{\bf 2}. Multiply each row containing at least three $-1$s  by $-1$. It follows that  there are at most two $-1$s in each row of $A$. 

{\bf 3}. If there is a row in $A$ with two $-1$s we permute it with the second row. The obtained matrix is permutationally equivalent to~$B$. 

{\bf 4}. Otherwise each row of $A$ contains at most one $-1$. $A$ is nonsingular. Hence there is no equal rows in $A$, and thus by Item {\bf 1} each row from second till $5$th contains $-1$. Then $\nu(A)=4$. By the nonsingularity of $A$, all $-1$s  are located in different columns. Hence by Lemma \ref{lem-1-1} $A$ is permutationally equivalent to $D_{(5,4)}$.
\end{prf}

\subsection{Matrices of size 5 and 6}

For $n = 5$ and $n = 6$ there are matrices which cannot be reduced to $D_{(n,n-1)}$, $D_{(n,n)}$ or a matrix  satisfying the condition ${\mathfrak A}$. Thus for nonsingular matrices we consider these two cases separately in this subsection. 

\begin{lemma}
Let $A \in M_5(\pm 1)$ be a nonsingular matrix. Then $|\per A| \le \per D_{(5,4)}$. Moreover, if $|\per A| = \per D_{(5,4)}$, then $A$ can be reduced to $D_{(5,4)}$ by standard transformations. 
\label{lemma5.0}
\end{lemma} 

\begin{prf}
If $A$ can be reduced to $D_{(5,4)}$ by standard transformations, then $|\per A| = \per D_{(5,4)}$ and the lemma is proved. 

Suppose that $A$ cannot be reduced to $D_{(5,4)}$ by standard transformations. Then by Lemma \ref{d5}    $A(3,4,5|)$ is reduced to
$ \left(
\begin{smallmatrix}
1 & 1 & 1 & 1 & 1 \\
1 & 1 & 1 & -1 & -1
\end{smallmatrix} \right).
$

Let us apply the generalized Laplace formula to the matrix $A$:
\begin{equation}
\per A = \suml{\al \in \SS_{5,2}}{} \per A[1,2|\al] \: \per A(1,2|\al).
\label{eq-p5}
\end{equation}

Since $\per \left(\begin{smallmatrix} 1 & 1 \\ 1 & -1 \end{smallmatrix} \right) = 0$, there are at most four non-zero summands in \eqref{eq-p5} which correspond to the submatrices 
\begin{equation}
A(1,2|1, 2), \ A(1,2|1, 3),\ A(1,2|2, 3),\ A(1,2|4, 5).
\label{eq-p5.2}
\end{equation}

Three possible options appear. 

{\bf 1.} There are no matrices of rank $1$ among the matrices listed in \eqref{eq-p5.2}. Then by Theorem~\ref{th-kr-new}  we obtain $|\per A(1,2|i, j)| < 6$ for all pairs $(i,j)\in\{(1,2),\: (1,3),\: (2,3),\: (4,5)\}$. Since for any matrix $B \in M_3(\pm 1)$ we have $\per B \in \{\pm 2, \pm 6\}$, the following inequality is obtained:
$
|\per A| \le 
2(|\per A(1,2|1, 2)| + |\per A(1,2|1, 3)| + |\per A(1,2|2, 3)| + |\per A(1,2|4, 5)|)  
= 16 < 24 = \per D_{(5,4)}.
$

{\bf 2. } One of the submatrices $A(1,2|1, 2), \ A(1,2|1, 3),\ A(1,2|2, 3)$ is of the rank $1$. Without loss of generality  $\rk A(1,2|1, 2) = 1$. It follows that columns of $A(1,2|1, 2)$ are proportional. 
Since $A$ is nonsingular,   $\rk A(1,2|) = 3$ and $\rk A(1,2|4,5) = 3$. Hence, $\rk A(1,2|1, 3) \ge 2$ and $\rk A(1,2|2, 3) \ge 2$.

We can multiply the rows of $A$ by $-1$ in such a way that the third column of $A$ does not contain $-1$s. Since $\rk A(1,2|1, 2) = 1$, we get that the entries in the second column of $ A(1,2|1, 2)$ are all the same and the entries in the third column of $ A(1,2|1, 2)$ are all the same. Then up to the permutation of the fourth and fifth columns of $A$ we have that  $A$ is of the following form:
 $$
A_1 = 
\left (\begin{smallmatrix}
1 & 1 & 1 & 1 & 1 \\
1 & 1 & 1 & -1 & -1 \\
a_{3,1} & a_{3,2} & 1 & 1 & 1 \\
a_{4,1} & a_{4,2} & 1 & 1 & 1 \\
a_{5,1} & a_{5,2} & 1 & 1 & 1 \\
\end{smallmatrix} \right),    
A_2 = 
\left (\begin{smallmatrix}
1 & 1 & 1 & 1 & 1 \\
1 & 1 & 1 & -1 & -1 \\
a_{3,1} & a_{3,2} & 1 & 1 & -1 \\
a_{4,1} & a_{4,2} & 1 & 1 & -1 \\
a_{5,1} & a_{5,2} & 1 & 1 & -1 \\
\end{smallmatrix} \right),     
A_3 = 
\left (\begin{smallmatrix}
1 & 1 & 1 & 1 & 1 \\
1 & 1 & 1 & -1 & -1 \\
a_{3,1} & a_{3,2} & 1 & -1 & -1 \\
a_{4,1} & a_{4,2} & 1 & -1 & -1 \\
a_{5,1} & a_{5,2} & 1 & -1 & -1 \\
\end{smallmatrix} \right),
$$
where  $a_{i,j} \in \{\pm 1\}$ and $\rk A(1,2|4,5)=3$. 
The fourth and fifth columns of $A_1$ and $A_3$ are equal which contradicts to the invertibility of $A$. Hence~$A = A_2$. 

We multiply the last column of $A$ by $-1$. Let us consider the pairs $(a_{i,1}, a_{i,2}) \in \{\pm 1\}^2$, where $i \in \{3,4,5\}$. There are four such pairs. Since all rows of $A(1,2|1, 2)$ are the same and $A$ is invertible, it follows that all rows of $A(1,2|3,4,5)$ are distinct. So, we have to choose any three different pairs from the given pairs. Since we do not distinguish permutationally equivalent matrices, only the following   variants  occur:
$
A_{2,1} = 
\left (\begin{smallmatrix}
1 & 1 & 1 & 1 & -1 \\
1 & 1 & 1 & -1 & 1 \\
1 & -1 & 1 & 1 & 1 \\
-1 & 1 & 1 & 1 & 1 \\
1 & 1 & 1 & 1 & 1 \\
\end{smallmatrix} \right),$  $  
A_{2,2} = 
\left (\begin{smallmatrix}
1 & 1 & 1 & 1 & -1 \\
1 & 1 & 1 & -1 & 1 \\
1 & -1 & 1 & 1 & 1 \\
-1 & 1 & 1 & 1 & 1 \\
-1 & -1 & 1 & 1 & 1 \\
\end{smallmatrix} \right),  
$ 
$
A_{2,3} = 
\left (\begin{smallmatrix}
1 & 1 & 1 & 1 & -1 \\
1 & 1 & 1 & -1 & 1 \\
1 & -1 & 1 & 1 & 1 \\
1 & 1 & 1 & 1 & 1 \\
-1 & -1 & 1 & 1 & 1 \\
\end{smallmatrix} \right), $ or $
A_{2,4} = 
\left (\begin{smallmatrix}
1 & 1 & 1 & 1 & -1 \\
1 & 1 & 1 & -1 & 1 \\
-1 & 1 & 1 & 1 & 1 \\
1 & 1 & 1 & 1 & 1 \\
-1 & -1 & 1 & 1 & 1 \\
\end{smallmatrix} \right).
$
$A_{2,1}$ is permutationally equivalent to  $D_{(5,4)}$ and it is straightforward to check that  $\per A_{2,2} = \per A_{2,3} = \per A_{2,4} = 8 < \per D_{(5,4)}.$
Thus in this case the maximum of $|\per A|$ can be achieved only on a matrix which can be reduced to $D_{(5,4)}$ by standard transformations.

{\bf 3.}  $\rk A(1,2|4,5) = 1$. Then the columns of $A(1,2|4,5)$ are proportional. We multiply the rows $3, 4, 5$ of $A$ by either 1 or $-1$ in order to obtain that $a_{3,1}=a_{4,1}=a_{5,1}=1$. Hence, for each column of $A(1,2|4,5)$ it holds that all its entries are either all positive or all negative. Then by the pigeonhole principle at least two columns are equal.  Since the submatrix $A[1,2|1,2,3]$ contains only $1$, at least two of the first three columns of $A$ are equal. Therefore,  $A$ is singular. This contradiction  concludes the proof. 
\end{prf}

\begin{lemma}
Let $A \in M_6(\pm 1)$ be nonsingular matrix, and assume that Conjecture \ref{h1} holds for all matrices of orders less than 6. 
 Then $|\per A| \le \per D_{(6,5)}$. Moreover, if $|\per A| = \per D_{(6,5)}$, then $A$ can be reduced to $D_{(6,5)}$ by standard transformations.
\label{lemma6.0}
\end{lemma} 

\begin{prf}
By Lemma \ref{lem:vid} the matrix $A$ can be reduced by standard transformations to one of the following forms:

{\bf 1.} The matrix $D_{(6,5)}$. In this case $|\per A| = \per D_{(6,5)}$. 

{\bf 2.} The matrix $D_{(6,6)}$. In this case by Lemma \ref{PropSmallD} $|\per A| = \per D_{(6,6)} = 112 < 128 = \per D_{(6,5)}$. 

{\bf 3.} The matrix $P_1$ or   $P_2$. In both cases  $|\per A| = 16 < \per D_{(6,5)}$. 

{\bf 4.} Matrix $A$ satisfies the condition ${\mathfrak A}$. Hence, without loss of generality we obtain:
$
A(3,4,5,6|) = 
\left(\begin{smallmatrix}
1 & 1 & 1 & 1 & 1 & 1 \\
1 & 1 & 1 & -1 & -1 & -1
\end{smallmatrix}
\right)$.

Let us apply the generalized Laplace formula to the matrix $A$:
\begin{equation}
\per A = \suml{\al \in \SS_{6,2}}{} \per A[1,2|\al] \per A(1,2|\al)
\label{eq-p6}
\end{equation}

Since $\per\left(\begin{smallmatrix} 1 & 1 \\ -1 & 1\end{smallmatrix}\right) = 0$, we obtain that there are at most six non-zero summands in the expression~\eqref{eq-p6}. 

The matrix $A$ is nonsingular. Hence $\rk A(1,2|)=4$. Then there is a $4\times 4$ submatrix in $A(1,2|)$ of the full rank and there is no $4\times 4$ submatrix in $A(1,2|)$ of the rank $1$. Let the indices $i,j,k,l$ be such that $\rk A[3,4,5,6|i,j,k,l] = 4$. Since  there are three $1$s and three $-1$s in the second row of $A$, we get that there are both $1$s and $-1$s  in set $\{a_{2,i}, a_{2,j}, a_{2,k}, a_{2,l}\}$.  Up to the multiplication of the second row by $-1$ two options appear:

{\bf 4.1.} There are three   $1$s among the entries $a_{2,i}, a_{2,j}, a_{2,k}, a_{2,l}$. Up to the permutation of the columns 4, 5 and 6 we assume that $\rk A(1,2|5,6) = 4$. Then for $A'\in \{A(1,2|4,5), \; A(1,2|4,6)\}$ it holds that $\rk A'\ge 3$.  Since Conjecture \ref{h1} holds for singular matrices of order $4$, we obtain that if $\rk A'=3$,  then $|\per A'|\le \per D_{(4,2)}=8$ by Proposition~\ref{PropSmallD}. If $\rk A'=4$,  then by Theorem \ref{experiment} we obtain that $|\per A'|\le 8$ and $|\per A(1,2|5,6)| \le 8$. Summarizing, we have $|\per A(1,2|\al_1, \al_2)| \le 8$ where $\al_1,\al_2 \in \{4,5,6\}$. 

The rest $4\times 4$ submatrices have the form $A''= A(1,2|\be_1, \be_2)$, where $\be_1, \be_2 \in \{1,2,3\}$. By the aforesaid $\rk A''\ge 2$. Since Conjecture \ref{h1} holds for singular matrices of order 4, we obtain $|\per A''|\le \per D_{(4,1)}=12$  by Proposition \ref{PropSmallD} if $A''$ is singular. If $A''$ is nonsingular, then by Theorem \ref{experiment} we obtain $|\per A''| \le 8< 12$, so in both cases $|\per A''|\le 12$ holds. Now using the equality \eqref{eq-p6} we obtain that
$
|\per A| \le 2(3\cdot 8 + 3\cdot 12) = 120 < 128 = \per D_{(6,5)}. 
$

{\bf 4.2.} There are two $1$s and two  $-1$s among the entries $a_{2,i}, a_{2,j}, a_{2,k}, a_{2,l}$. Up to the permutations of columns  we have $\rk A(1,2|1,6) = 4$. For the matrix $A''' \in \{A(1,2|1,2), A(1,2|1,3), A(1,2|4,6), A(1,2|5,6)\}$ it follows that $\rk A'''\ge 3$. Thus repeating the arguments of Item ${\bf 4.1}$  we get that $|\per A'''|\le 8$. Similarly for $A^*\in \{A(1,2|2,3),A(1,2|4,5)\}$ we have $\rk  A^*\ge 2$ and thus $|\per A^*|\le 12$. Now we use the equality \eqref{eq-p6} to obtain:
$
|\per A| \le 2(4\cdot 8 + 2\cdot 12) = 112 < 128 = \per D_{(6,5)}. 
$
The lemma is proved. 
\end{prf}

\subsection{Matrices of size greater than 6}

Here we prove the induction step for nonsingular matrices of order $n \ge 7$. 

To prove the induction step for nonsingular matrices we need to find  the decomposition for the permanent of $D_{(n,n-1)}$ via the permanents of $D_{(k,m)}$ for smaller $k$ and $m$.

\begin{lemma} \label{lem:d}  Let $D=D_{(n,n-1)} \in M_n(\pm 1)$ and $n \ge 5$. Then 
\begin{equation}
\per D = 2 \per D_{(n-2,n-3)} +
\label{eqm6.3}
\end{equation}
$$
 (n^2-7n+12)\per D_{(n-2,n-5)} + 2(n-3)\per D_{(n-2,n-4)}.
$$
\end{lemma}

\begin{prf}
We apply the generalized Laplace formula to the matrix  $D = D_{(n,n-1)}$.
 The first two rows of this matrix look as 
$
\left(\begin{smallmatrix}
	-1 & 1 & 1 & \ldots & 1 \\
	1 & -1 & 1 & \ldots & 1 
\end{smallmatrix}\right).
$ 

All  submatrices of size $2$ of the above matrix can be considered as follows:

{\bf 1.} The submatrix  $D[1,2|1,2]$. Its permanent is equal to  $2$. The complementary submatrix is equal to $D_{(n-2, n-3)}$.
	
{\bf 2.} Submatrices  $D[1,2|i,j]$ where $i \in \{1,2\}, j \in \{3, \ldots, n\}.$ The number of such submatrices is equal to $2(n-2)$. The permanent of any of these matrices is equal to zero. 
	
{\bf 3.} Submatrices $D[1,2|i,j]$ where both $i,j \in \{3, \ldots, n-1\}$. The number of such submatrices is $\frac12(n-3)(n-4) = \frac12(n^2 - 7n + 12)$. The permanent of each  of these submatrices is equal to  $2$. Moreover in the complementary submatrix $D(1,2|i,j)$ we have  $n-5$ negative entries with no two of them lying in the same line (row or  column). 
By Lemma \ref{lem-1-1} the matrix $D(1,2|i,j)$ is permutationally equivalent to  $D_{(n-2, n-5)}$ and we obtain $\per D(1,2|i,j) = \per D_{(n-2, n-5)}$. 
	
{\bf 4.} Submatrices  $D[1,2|i,n]$, where  $i \in \{3, \ldots, n-1\}$. There are  $n-3$ such submatrices. Permanent of each of them is equal to $2$. Moreover in the complementary submatrix  $D(1,2|i,n)$ we have  $n-4$ negative entries with no two of them are in the same line. 
By Lemma \ref{lem-1-1} the matrix $D(1,2|i,n)$ is permutationally equivalent to  $D_{(n-2, n-4)}$ and $\per D(1,2|i,n) = \per D_{(n-2, n-4)}$. 

Now considering separately the above groups of summands in the generalized Laplace expansion formula
 we obtain: 
$$
\per D = \suml{\al \in \SS_{n,2}}{} \per D[1,2|\al] \per D(1,2|\al) = \\
$$
$
 = 2\per D_{(n-2,n-3)} + (n^2-7n+12)\per D_{(n-2,n-5)} + 2(n-3)\per D_{(n-2,n-4)}.
$
\end{prf}

\begin{theorem}
Let $A \in M_n(\pm 1)$ be nonsingular, $n \ge 7$. Assume that Conjecture \ref{h1} is true for all
matrices of the order $m<n$. 
Then the inequality $|\per A| \le \per D_{(n,n-1)}$ holds. Moreover, if $|\per A| = \per D_{(n,n-1)}$, then $A$ can be reduced to $D_{(n,n-1)}$ by standard transformations.  
\label{theorem1.0}\end{theorem}

\begin{prf}
Applying Lemma \ref{lem:vid} to the matrix  $A$ we get just 3 types of matrices to be considered, since $P_1,P_2 \in M_6(\pm 1)$, and $n \ge 7$ in this lemma. Let us investigate these types separately. 

{\bf 1.} Assume that $A$ can be reduced to $D_{(n,n-1)}$ by standard transformations. Then   $|\per A| = \per D_{(n,n-1)}$ since standard transformations preserve module   of the permanent. It remains to show  that in other cases  $\per A < \per D_{(n,n-1)}$. 

{\bf 2.} Suppose, $A$ can be reduced to the form  $D_{(n,n)}$. Then we have the equality  $|\per A| = \per D_{(n,n)}$. By Lemma  \ref{lem40} we have the inequality  $\per D_{(n,n)} < \per D_{(n,n-1)}$. Thus $|\per A| < \per D_{(n,n-1)}$, as required.

{\bf 3.} Finally, consider the case that the condition ${\mathfrak A}$ is satisfied for $A$. It means that   $A$ has the following form:

 $\bullet$ The first row consists of $1$s. 

 $\bullet$ The second row has at least three $1$s and three $-1$s. 
 
Assume that in the second row there are  $k$ negative entries, $3\le k\le n-3$. Let us apply the generalized Laplace formula to the matrix $A$:
\begin{equation}
\per A = \suml{\al \in \SS_{n,2}}{} \per A[1,2|\al] \per A(1,2|\al)
\label{eqm6}
\end{equation}

Since  $A$ is nonsingular, $\rk A(1,2|\al)\ge n - 4$.
Since $n \ge 7$, the inequality in Conjecture \ref{h1} holds for matrices of the order  $n-2 > 4$. Hence by Corollary \ref{con40} there is a bound:  
\begin{equation}
|\per A(1,2|\al)| \le \per D_{(n-2,n-5)}.
\label{eqm6.1}
\end{equation}
Now we estimate the permanent of  $A$. Consider all submatrices $A[1,2|\al]$ of order $2$ where $\al \in \SS_{n,2}$. For each of them we have that both entries in the first row are 1.
For the second row we have the following alternatives:

$\bullet$ There are two  $1$s in the second row. Then $\per A[1,2|\al] = 2$.

$\bullet$ There are two  $-1$s   in the second row. Then $\per A[1,2|\al] = -2$.

$\bullet$ The second row contains the entries $1$ and $-1$. Then $\per A[1,2|\al] = 0$.

If $\per A[1,2|\al] = 0$, then the corresponding summand in  \eqref{eqm6} is zero. The number of zero summands in  \eqref{eqm6} is equal to  $C_{n-k}^1C_k^1 = k(n-k)$ where $k$ is number of $-1$s in the second row. Taking into account all the possible values of the parameter $k$ we obtain that the number of 0 summands in \eqref{eqm6} is minimal if either $k = 3$ or $k = n-3$. In both cases there are $3(n-3)$ summands which are equal to~0.

By \eqref{eqm6.1} each non-zero summand in  \eqref{eqm6} can be estimated by  $2\,\per D_{(n-2,n-5)}$. 
Altogether there are $C_n^2 = \frac12n(n-1)$ summands in \eqref{eqm6}. 
Thus the permanent function of the matrix $A$ can be estimated by
\begin{equation}
|\per A| \le (\frac12n(n-1) - 3(n-3)) \cdot 2\cdot \per D_{(n-2,n-5)} = 
\label{eqm6.2}
\end{equation}
$$
(n^2 - 7n + 18)\per D_{(n-2,n-5)}.
$$

Let us estimate the difference between the permanent of  $D_{(n,n-1)}$ and the permanent of $A$. Using the exact formula \eqref{eqm6.3} computed in Lemma~\ref{lem:d} and the bound from the formula \eqref{eqm6.2} we obtain that
$$\per D_{(n,n-1)} - |\per A| \ge $$
\begin{equation}
\ge 2\per D_{(n-2,n-3)} + 2(n-3)\per D_{(n-2,n-4)} - 6 \per D_{(n-2, n-5)}. 
\label{eqm6.4}
\end{equation}

By Lemma \ref{lem:c} we have the equality: 
\begin{equation}
\per D_{(n-2,n-5)} = \per D_{(n-2,n-4)} + 2 \per D_{(n-3,n-5)}.
\label{eqm6.5}
\end{equation}

Substitute  \eqref{eqm6.5} into  \eqref{eqm6.4} and obtain:
$$
\per D_{(n,n-1)} - |\per A| \ge
$$
\begin{equation}
\ge 2\,\per D_{(n-2,n-3)} + 2(n-6)\per D_{(n-2,n-4)} - 12\,\per D_{(n-3, n-5)}. 
\label{eqm6.6}
\end{equation}

Applying to  $\per D_{(n-2,n-3)}$ the result of Lemma  \ref{lem:c} we obtain 
$$
\per D_{(n-2,n-3)} = \per D_{(n-2,n-2)} + 2\per D_{(n-3,n-3)} . 
$$
Now applying Lemma  \ref{lem:c2} to the first summand we obtain
$$
\per D_{(n-2,n-3)}= (n-2)\per D_{(n-3,n-3)} + 2(n-3)\per D_{(n-4,n-4)} . 
$$
Then by  Lemma \ref{lem:c2} applied to the first summand of the obtained formula we have 
\begin{equation}
\per D_{(n-2,n-3)}= 
\label{eqm6.7}
\end{equation}
$$
(n^2 - 5n + 4)\per D_{(n-4,n-4)} + 2(n^2 - 6n + 8)\per D_{(n-5,n-5)}.
$$

Now let us apply Lemma  \ref{lem:c} to   $\per D_{(n-3,n-5)}$. We get
$$
\per D_{(n-3,n-5)} = \per D_{(n-3, n-4)} + 2\per D_{(n-4,n-5)} . 
$$
Applying Lemma  \ref{lem:c} to  each summand we obtain
$$
\per D_{(n-3,n-5)} =
\per D_{(n-3,n-3)} + 4\per D_{(n-4,n-4)} + 4\per D_{(n-5,n-5)} .
$$
Then applying Lemma  \ref{lem:c2} to the first  summand we get the equality: 
\begin{equation}
\per D_{(n-3, n-5)} = (n-1)\per D_{(n-4,n-4)} + (2n - 4)\per D_{(n-5,n-5)}.
\label{eqm6.8}
\end{equation}

By  Lemma \ref{lem:c} 
$
\per D_{(n-2,n-4)} = \per D_{(n-2, n-3)} + 2 \per D_{(n-3, n-4)} = 
$
$
\per D_{(n-2, n-2)} + 4\per D_{(n-3, n-3)} + 4 \per D_{(n-4,n-4)},  
$
for the second equality we applied  Lemma \ref{lem:c}  to each summand and collected common summands.
Now apply Lemma  \ref{lem:c2} to the first summand and recollect the summands to obtain
$
\per D_{(n-2,n-4)}= n\per D_{(n-3,n-3)} + (2n - 2)\per D_{(n-4,n-4)}.
$
Applying Lemma  \ref{lem:c2} to the first summand and reordering we get
\begin{equation}
\per D_{(n-2,n-4)} = (n^2-3n-2)\per D_{(n-4,n-4)} + (2n^2 - 8n)\per D_{(n-5,n-5)}
\label{eqm6.9}
\end{equation}

By substitution of the expressions  \eqref{eqm6.7}, \eqref{eqm6.8} and \eqref{eqm6.9} into  the inequality \eqref{eqm6.6} we  obtain the following inequality:
$$
\per D_{(n,n-1)} - |\per A| \ge (2n^3 - 16n^2 + 10n + 44)\per D_{(n-4,n-4)} +
$$
\begin{equation}
 + (4n^3 - 36n^2 + 48n + 80) \per D_{(n-5,n-5)} = F(n).
\label{eqm6.10}
\end{equation}

Consider the polynomial function $f(n) = 2n^3 - 16n^2 + 10n + 44$. Its derivative $f'(n) = 6n^2 - 32n + 10$ has the roots $\frac{1}{3}$ and $5$. Since $f(7) = 16$, it follows that  $f(n) > 0$ for any $n \ge 7$. 

Consider also the polynomial function  $g(n) = 4n^3 - 36n^2 + 48n + 80$. Its derivative $g'(n) = 12n^2 - 72n + 48$ has the roots $3 \pm \sqrt{5}$. Since $g(7) = 24$ and $3 \pm \sqrt{5}<7$, it follows that  $g(n) > 0$ for any $n \ge 7$. 

By Proposition~\ref{PropSmallD} for $n \ge 9$ the inequalities $\per D_{(n-4,n-4)} > 0$ and $\per D_{(n-5,n-5)} > 0$ hold. Therefore, $F(n) > 0$ for any $n \ge 9$. 

By direct computation we obtain $F(8) = 576>0$, $F(7) = 17>0$. Hence, $F(n) > 0$ for any $n \ge 7$. 

Since $F(n) > 0$ for $n \ge 7$, the obtained bound is strict and $A$ cannot be reduced to $D_{(n,n-1)}$.
\end{prf}

Thus we proved the following theorem:

\begin{theorem}
Let $A \in M_n(\pm 1)$ be nonsingular and $n \ge 5$. Assume that Conjecture \ref{h1} is true for all
matrices of the order $m$ for all $m<n$.  
Then the inequality $|\per A| \le \per D_{(n,n-1)}$ holds. Moreover if $|\per A| = \per D_{(n,n-1)}$, then  $A$ can be reduced to $D_{(n,n-1)}$ by standard transformations.
\label{theorem1}\end{theorem}

\begin{prf}
The result follows  for $n=5$ from Lemma \ref{lemma5.0},  for $n=6$ from Lemma \ref{lemma6.0}, and for $n\ge 7$ from Theorem \ref{theorem1.0}.
\end{prf}

\section{Properties of the rank vectors} \label{S5}

In the following lemmas we use the notations introduced in Section 2, see Definitions \ref{def_rk_vect} --- \ref{defCb}. Everywhere in this and the next  sections we consider sets which may contain repeating elements. Also considering unions of such sets we say that they also contain repeating elements and the number of such elements is equal to the sum of the numbers of  corresponding elements in each of the sets. For example, $\{a\}\cup\{a\}=\{a,a\}$. 

\begin{lemma}
 Let $Y, Z\subseteq M_k(\pm 1)$ be matrix subsets. We assume that repetitions of matrices in these sets are possible. Let $X = Y \cup Z$. Then the equality $R(X) = R(Y) + R(Z)$ holds true, where $R(X)$ is introduced in Definition~\ref{def_rk_vect}. 
\label{lem:setsum}\end{lemma}

\begin{prf}
The sum $R(Y) + R(Z)$ is defined correctly since the vectors  $R(Y)$  and $R(Z)$ 
are of the same size. By definition, the  $i$th component of the vector $R(X)$ is equal to the number of matrices  of the rank  $k-i+1$. Since $X$ is the union of  $Y$ and $Z$, the number of matrices of rank  $k-i+1$  in $X$ is equal to the sum of the numbers of matrices with the same rank  in $Y$  and $Z$. Therefore the equality  $R(X) = R(Y) + R(Z)$ is true   for the $i$th coordinate and any $1 \le i \le k$  which proves the statement.
\end{prf}

\begin{lemma}
Let $A \in M_{k,n}(\pm 1)$, $b\in M_{k,1}(\pm 1)$ and $B\in M_{k,n+1}(\pm 1)$ be the  concatenation  of $A$ with the  column  $b$. Then   $K(A) \subset K(B)$ holds.
\label{con:setin}\end{lemma}
\begin{prf}
By definition   the set $K(A)$  contains all submatrices of  $A$ of the order $k$. Since  $A$ is a submatrix of  $B$ and $K(B)$ consists of all $k \times k$ submatrices  of $B$, we have $K(A) \subset K(B)$. 
\end{prf}

\begin{lemma}
Let $A \in M_{k,n}(\pm 1)$, $b\in M_{k,1}(\pm 1)$ and $B\in M_{k,n+1}(\pm 1)$ be the  concatenation  of $A$ with the  column  $b$.  Then  the following inclusion holds:
\begin{equation}
K(B) \subset K(A) \cup (\bigcup_{C \in K(A)} <C,b>).
\label{eq:vkl}
\end{equation}
\end{lemma}

\begin{prf}
By Lemma~\ref{con:setin}, $K(A) \subset K(B)$. To show the inclusion \eqref{eq:vkl}, it is sufficient to prove that for any matrix $X \in (K(B)\setminus K(A))$ there is $C \in K(A)$ such that $X \in <C,b>$.

Let  $C' \in (K(B)\setminus K(A))$. Since $C'$ is not a submatrix of $A$, it contains the column $b$. All other columns of $C'$ are columns of $A$. It follows that there is matrix $C \in K(A)$ such that up to a permutation of columns it differs from $C'$ by one column which corresponds to $b$. Thus, we have $C' \in <C, b>$.  Since $C'$  is chosen arbitrarily, the inclusion \eqref{eq:vkl}  is proved. 
\end{prf}

\begin{lemma}
Let $A \in M_{k,n}(\pm 1)$, $b\in M_{k,1}(\pm 1)$ and $B\in M_{k,n+1}(\pm 1)$ be the  concatenation  of $A$ with the  column  $b$. Then any matrix from $K(B) \setminus K(A)$  occurs exactly  $(n-k+1)$ times in the set  $\cup_{C \in K(A)} <C,b>$. Moreover, the set  $\cup_{C \in K(A)} <C,b>$ does not contain other matrices. 
\label{lem:sub}
\end{lemma}

\begin{prf}
Any matrix $C' \in (K(B)\setminus K(A))$ contains the vector $b$. All other $k-1$ columns of the matrix $C'$ are columns of $A$. Denote the matrix obtained by these $k-1$ columns by $C$. Therefore, if we take any matrix from 
 $K(A)$, such that the columns of $C$ are fixed and one column $c$ is chosen randomly, we can obtain  $C'$ by replacing the column $c$ by the column $b$. Totally, there are $n-k+1$ such matrices  $C \in K(A)$. Therefore, any matrix from  $K(B) \setminus K(A)$ occurs exactly $n-k+1$ times in the set  $\cup_{C \in K(A)} <C,b>$.

It remains to note that any matrix from $<C, b>$ lies in $K(B)\setminus K(A)$ since it is a submatrix in $B$, which contains 
the column $b$. This proves the lemma.
\end{prf}

\begin{cons}
Let $A \in M_{k,n}(\pm 1)$, $b\in M_{k,1}(\pm 1)$ and $B\in M_{k,n+1}(\pm 1)$ be the  concatenation  of $A$ with the  column  $b$. Then we have: 
$$
R(K(B)\setminus K(A)) = \frac{1}{n-k+1}R(\bigcup_{C \in K(A)} <C,b>)
$$
\label{con:veq}
\end{cons}

\begin{prf}
By Lemma  \ref{lem:sub} the set  $\cup_{C \in K(A)} <C,b>$ can be represented  as the union of  $n-k+1$ copies of the set $K(B)\setminus K(A)$. The application of Lemma  \ref{lem:setsum} leads to the equality 
$$
(n-k+1)R(K(B)\setminus K(A)) = R(\bigcup_{C \in K(A)} <C,b>).
$$

Dividing by $n-k+1$, we obtain the result. 
\end{prf}

Our next goal is to describe how the rank of the matrix  $C \in M_k(\pm 1)$ is connected to the ranks of matrices from 
the set $<C, b>$ where the vector $b$ is of the dimension $k$ and consists of signed ones.
We have the following characterization. 

\begin{lemma}
Let $C \in M_k(\pm 1),$   $\rk C = m \ge 1$, and $b \in M_{k,1}(\pm 1)$. Then  
one of the following statements holds. 

{\bf 1.} If the column  $b$ is a linear combination of  $l \le m$ columns of the matrix  $C$ and  $l$ is minimal, then for any $A\in <C,b>$ either $\rk A =\rk C$ or $\rk A = \rk C - 1$ and the number of matrices $A \in <C,b>$ such that $\rk A = \rk C-1$ is less than or equal to $m-l$.

{\bf 2.} If the column  $b$ and the columns of the matrix $C$ are linear independent, then for any $A\in {<C,b>}$ either $\rk A =\rk C$ or  $\rk A=\rk C+1$.
\label{lem:rank}\end{lemma}

\begin{prf}
Since rank of a matrix does not depend on the order of its columns, the set of ranks of matrices from the set  $<C,b>$ does not depend on the order of columns in $C$. The rank of matrix  $C$ is equal to $m$, therefore, we can reorder its columns in such a way that the first  $m$ columns are linear independent and the others are their linear combinations. In the sequel, we assume that this condition is automatically satisfied for the matrix $C$. 
 
{\bf 1. } Let $b$ be  a linear combination of the columns of  $C$. Without loss of generality we can assume that $b$ is a linear combination of the first  $l \le m$ columns of the matrix  $C$ where  $l$ is the minimal number of vectors
which are necessary to  represent the vector $b$. Indeed, by the minimality of  $l$ the columns of $C$ which represent vector $b$ are linear independent. Let us put these vectors on the first $l$ positions in matrix $C$ and 
complement them by $m-l$ linear independent columns of $C$.
We obtain that any of the remaining columns of  $C$ is a linear combination of the first  $m$ columns. 
There are several variants to be considered: 

{\bf A.} The vector $b$ replaces one of the first  $l$ columns of   $C$ and the obtained matrix is  $C_i$,
where $i \le l$. In this case the subspace spanned by the first  $l$ columns of    $C$ coincides with the space
spanned by the first  $l$ columns of $C_i$. Therefore,   $C_i$ can be reduced to  $C$ by means of the elementary transformations. It follows that $\rk C = \rk C_i$, where $i \le l$. 
	
{\bf B.} The vector  $b$ replaces $i$th column of $C$ where $i \in \{ l+1, \ldots, m\}$. In this case the matrix 
$C_i$ is a concatenation of the matrix  $C$ and the vector $b$ and deleting one of the columns. Since $b$ is a linear combination of columns of $C$, we have that $rk(C\vert b)=\rk (C)$. Elimination of one of the columns can either leave the rank unchanged or decrease it by 1.  So, $\rk(C_i)\in\{\rk(C),\rk(C)-1\}$.
	
{\bf C.} The vector $b$ replaces $i$th column of $C$ where $m+1 \le i \le k$. 
The column which was replaced by $b$ was a linear combination of the first  $m$ columns.
The vector  $b$ can be expressed as a linear combination of the first  $m$ columns of $C$. Therefore, for the obtained matrix $C_i$ we have $\rk C_i = \rk C$.  

Thus the matrix  $C_i \in <C, b>$ of rank  $(\rk C - 1)$ can be obtained  only by the replacement 
of  one of the columns with indices  $l+1, \ldots, m$ by the column  $b$. Altogether, we have less than or equal to $m-l$ different $C_i$ and the first part of the lemma is proved. 

{\bf 2. } Let us assume now that the column  $b$ cannot be represented as a linear combination of the columns of  $C$. Every matrix $C_i$ from  $<C, b>$ is obtained as follows: we add one column which is equal to $b$ and
after we remove one column except the added one. If we add a column  $b$, then the rank of the obtained matrix  $C'$ 
becomes greater by one, since this column cannot be expressed as a linear combination of the other columns of $C$. Removing of the column either does not change the  rank or reduces it by $1$.
Therefore, for any matrix  $C_i \in <C, b>$ its rank is equal to $\rk C$ or  $(\rk C + 1)$. 
\end{prf}

\begin{cons}
Let $C \in M_k(\pm 1)$, $\rk C = m > 1$ and $b$ be a vector of signed ones.
Then the number of matrices of the rank  $\rk C - 1$ in the set  $<C,b>$ for an arbitrary  $b$ 
is less than or equal to $m-1$. This bound can be achieved  only if
there exists a column equal to  $\pm b$  among the columns of~$C$. 
\label{con1}\end{cons}

\begin{prf}
By Lemma \ref{lem:rank} the number of matrices  $C_i \in <C, b>$ with $\rk C_i = \rk C - 1$ is not greater than  $m - l \le m-1$ since $l \ge 1$ where $m - l$ is the number of the columns of $C$ such that $b$ lies in the span of these columns. It is clear that  $m - l$ is maximal for  $l = 1$. This implies that the column $b$  must be a scalar multiple of one of the columns of  $C$. Since we consider  $(-1,1)$-matrices,  there is a column equal to  $\pm b$  in~$C$. 
\end{prf}



\begin{lemma}
Let $D = D_{(k,l)} \in M_k(\pm 1)$ where $k > l$. Then $\rk D = l + 1$.
\label{lem:rkD}\end{lemma}

\begin{prf}
Since $D[1,\ldots,l+1|1,\ldots,l+1] = D_{(l+1,l)}$ and $D_{(l+1,l)}$ is nonsingular,   $\rk D \ge l+1$. Note that the rows $l+1,\ldots, n$ in $D$ are equal. It follows that the rows $l+2,\ldots, n$ can be represented as a linear combinations of the  rows $1, \ldots, l+1$ and $\rk D \le l+1$. If we combine the two inequalities, we obtain $\rk D = l+1$.
\end{prf}

\begin{lemma}
Let $D = D_{(k,l)} \in M_k(\pm 1)$, $k > l$, and  $b=(1,\ldots,1)^t$ be of the size  $k$. Then $\rk D = (l + 1)$ and  there are  $l$ matrices of the rank  $l$ in  $<D,b>$. All other  matrices in  $<D,b>$ have the rank $(l+1)$.
\label{lem:kr}\end{lemma}

\begin{prf}
By Lemma \ref{lem:rkD}, $\rk D = l+1$. 

If we replace one of the columns $1, \ldots, l$ of $D$   by the column of 1, then, up to the permutation of rows and columns, we obtain the matrix  $D_{(k,l-1)}$. By Lemma \ref{lem:rkD}, $\rk D_{(k,l-1)}=l$. 

If we replace one of the columns $l+1, \ldots, k$ of $D$ by $b$, then $D$ does not change as a matrix in $M_k(\pm 1)$. Therefore, 
its rank is $(l+1)$. 
\end{prf}

\section{Constructing partial orders} \label{S6}

In this section  we introduce a partial order relation (a majorization) on a set of vectors.  
We show that the rank vector (in the sense of Definition~\ref{def_rk_vect}) of the matrix $D_{(n,k,k-1)} \in M_{k,n}(\pm 1)$ is minimal with respect to this order among all rank-vectors $R(A)$ where $A \in M_{k,n}(\pm 1)$  and $\rk A = k$.

\begin{defn} \cite{MaOlAr11}
Let $R_1 = (r^1_i)$ and $R_2 = (r^2_i)$ be real (actually integer) vectors of the same size. We say that $R_1\lep R_2$ if all following  inequalities are true:
\begin{equation}
\suml{i=1}{j} r^1_i \le \suml{i=1}{j} r^2_i \;\; \text{ for } \;\; j = 1, \ldots, k. 
\label{eq:por}
\end{equation}
We say that $R_{1} <^* R_2$ if $R_1\lep R_2$ and at least one inequality in \eqref{eq:por} is strict.
\end{defn}

\begin{lemma}
The formula  \eqref{eq:por} defines the partial order on the set of vectors from ${\mathbb R}^m$ for any $m>1$.
\end{lemma}

We plan to apply this order to work with the rank vectors determined in Definition~\ref{def_rk_vect}.

Next lemma is a consequence of  Lemma \ref{lem:setsum}. 

\begin{lemma}
Let $X_i$, $i = 1, \ldots, 4$ be some sets consisting of matrices of order $k$. Assume that $R_1 \lep R_2$ and $R_3 \lep R_4$ where $R_i = R(X_i) = (r^i_j)$, $i = 1, \ldots, 4$ and $j = 1, \ldots, k$. Then $R(X_1 \cup X_3) = R_1 + R_3 \lep R_2 + R_4 = R(X_2 \cup X_4)$. 
\label{lem:n1}
\end{lemma}

\begin{prf}
Indeed, by Lemma \ref{lem:setsum} we have $R(X_1 \cup X_3) = R_1 + R_3$ and $R(X_2 \cup X_4) = R_2 + R_4$. Denote by  $R(\cdot)_i$  the $i$th coordinate of the vector $R(\cdot)$. This leads to the following inequality: 
$$
\suml{i = 1}{j} R(X_1 \cup X_3)_i = \suml{i = 1}{j}r^1_i + \suml{i = 1}{j}r^3_i \le
\suml{i = 1}{j}r^2_i + \suml{i = 1}{j}r^4_i = \suml{i = 1}{j} R(X_2 \cup X_4)_i.
$$
Since the last inequality is true for any $j = 1, \ldots, k$, our lemma is proved.
\end{prf}

We prove two technical lemmas.

\begin{lemma}
Let $C \in M_k(\pm 1)$, $D = D_{(k,l)} \in M_k(\pm 1)$, $k > l$ and $\rk C \ge \rk D$. Let $b_1  \in M_{k,1}(\pm 1)$ be arbitrary, and $b_2 = (1,\ldots, 1)^t\in M_{k,1}(\pm 1)$. Then $R(<D,b_2>) \lep R(<C,b_1>)$.
\label{lem:n2}
\end{lemma}

\begin{prf}
Denote $m=\rk C$ and the vectors  $R(<C, b_1>)=(r^1_1, \ldots, r^1_k)$,  $R(<D, b_2>)=(r^2_1, \ldots, r^2_k)$. Then by Lemma \ref{lem:rkD} the inequality $m\ge l+1$ holds. By Lemma \ref{lem:kr} the vector  $R(<D, b_2>)$ has two nonzero coordinates. Namely, $r^2_{k-l}=k-l$ corresponds to 
the number of matrices of rank $(l+1)$ in the set $<D, b_2>$, and $r^2_{k-l+1}=l$
corresponds to the number of matrices of rank $l$ in the  set $<D, b_2>$. For $i\notin \{k-l,k-l+1\} $ it holds that $r^2_i=0$. Similarly, by Lemma \ref{lem:rank} there are 3 nonzero coordinates in vector  $R(<C, b_1>)=(r^1_1, \ldots, r^1_k)$.  These three nonzero  coordinates   correspond to the numbers of matrices of the ranks  $(m+1)$, $m$, and $(m - 1)$ respectively. Moreover, the sum of the entries in each vector from $R(<D, b_2>)$ and $R(<C, b_1>)$ is equal to $k$. So, $r^1_j=0$ for $j\notin \{k-m,k-m+1,k-m+2\}$, $r^1_{k-m}+r^1_{k-m+1}+r^1_{k-m+2}=k$. The situation splits into the following 3 cases.

{\bf 1}. If $m > l+2$, then all positive entries of $R(<C, b_1>)$ have the indices less than the indices of the positive entries of $R(<D, b_2>)$. Thus, the inequalities \eqref{eq:por} for the partial sums of the entries of the vectors $R(<C, b_1>)$ and $R(<D, b_2>)$ hold true and $R(<D,b_2>) \lep R(<C,b_1>)$.

{\bf 2}. Suppose that $m = l+2$. Let us write the vectors $R(<C, b_1>) $ and $R(<D, b_2>) $ one over another and obtain the following matrix:
$$
\left(\begin{smallmatrix}
0 & \ldots & 0 & r_{k - l - 2}^1 & r_{k - l - 1}^1 & r_{k - l}^1 & 0 & 0 & \ldots & 0 \\
0 & \ldots & 0 & 0 & 0 & r_{k-l}^2 & r_{k-l+1}^2 & 0 & \ldots & 0
\end{smallmatrix}\right).
$$ 

Since  $r_{k - l - 2}^1 + r_{k - l - l}^1 + r_{k - l}^1 = k = r_{k-l}^2 + r_{k-l+1}^2$ and all coefficients are nonnegative integers, for any $i$ we obtain
$\sum_{j=1}^i r_j^1 \ge \sum_{j=1}^i r_j^2 .$  Thus   $R(<D,b_2>) \lep R(<C,b_1>)$.

{\bf 3}. Let $m = l+1$.  Similarly to the previous item, let us write vectors  $R(<C, b_1>) $ and $R(<D, b_2>)  $ into the following matrix:
\begin{equation} \left( 
\begin{smallmatrix}
0 & \ldots & 0 & r_{k - l - 1}^1 & r_{k - l}^1 & r_{k - l + 1}^1 & 0 & \ldots & 0 \\
0 & \ldots & 0 & 0 & r_{k-l}^2 & r_{k-l+1}^2 & 0 & \ldots & 0
\end{smallmatrix}\right).
\label{eq3006}\end{equation}
By Corollary \ref{con1} we have $r_{k - l + 1}^1 \le l=r_{k-l+1}^2 $. Then  $r_{k - l - 1}^1 + r_{k - l}^1 \ge r_{k-l}^2$. Therefore, $\sum_{j=1}^{k-l} r_j^1 \ge \sum_{j=1}^{k-l} r_j^2 .$ The inequalities 
$\sum_{j=1}^i r_j^1 \ge \sum_{j=1}^i r_j^2 $
for $i=1,\ldots,  (k-l-1)$ are trivial, and for $i=k-l+1,\ldots, k$ one easily gets $\sum_{j=1}^i r_j^1 = \sum_{j=1}^i r_j^2 =k.$ Hence, $R(<D,b_2>) \lep R(<C,b_1>)$ and the lemma is proved in each of the cases. 
\end{prf}

\begin{cons}
Let $C \in M_k(\pm 1)$, $D = D_{(k,l)} \in M_k(\pm 1)$, $k > l$ and $\rk C \ge \rk D$. Let $b_1, b_2 \in M_{k,1}(\pm 1)$ and $b_2 = (1,\ldots, 1)^t$. 
Let $R(<C,b_1>) = (r_i)$ and $(<D,b_2>) = (u_i)$. Let $u_m$ be the last positive number in vector $(u_i)$. Then $u_m \ge r_m$ and if $u_m = r_m$, then 
 $R(<D,b_2>) = R(<C,b_1>)$.
\label{new-cons-2}\end{cons}

\begin{prf}
Since $C$ and $D_{(n,k,k-1)}$ are matrices of the same size, we obtain $\suml{i = 1}{k} r_i = \suml{i = 1}{k} u_i = S$. 
By Lemma \ref{lem:vec0} $R(C) \gep R(D_{(n,k,k-1)})$. Thus $\suml{i = 1}{k-1} r_i \ge \suml{i = 1}{k-1} u_i$ and we obtain the inequality
$
r_k = S - \suml{i = 1}{k-1} r_i \le S - \suml{i = 1}{k-1} u_i = u_k.
$

Suppose $r_m = u_m$. This equality possible only in Case {\bf 3} of the proof of the Lemma \ref{lem:n2}.
Considering \eqref{eq3006}, it follows that $r_{k-l+1}^2 = u_m = r_m = r_{k-l+1}^1$. 
By Corollary \ref{con1} the equality $r_{k-l+1}^2 = r_{k-l+1}^1$ is possible if and only if $b_1$ is proportional to one of the columns in matrix $C$. Thus in the set $<C,b_1>$ there is no matrices with rank greater than $\rk C$. It follows that $r_{k - l - 1}^1 = 0$ and $R(<D,b_2>) = R(<C,b_1>)$.
\end{prf}

\begin{lemma}
Let $A \in M_{k,n}(\pm 1)$,  $b\in M_{k,1}(\pm 1)$, and $B\in M_{k,n+1}(\pm 1)$ be the  concatenation  of $A$ with the  column  $b$. Assume that $k \le n$, $\rk A = \rk B = k$ and  $R(D_{(n,k,k-1)}) \lep R(A)$. Then $R(D_{(n+1,k,k-1)}) \lep R(B)$. Moreover, if $R(D_{(n,k,k-1)}) <^* R(A)$, then $R(D_{(n+1,k,k-1)}) <^* R(B)$.
\label{lem:vec0}\end{lemma}

\begin{prf}
By Lemma \ref{con:setin} we have the inclusions  $K(A) \subset K(B)$ and $K(D_{(n, k)}) \subset K(D_{(n+1, k)})$.

By Lemma  \ref{lem:setsum} the following equalities are true:  
\begin{equation}
R(B) = R\Bigl(K(A) \cup \bigl(K(B)\setminus K(A)\bigr)\Bigr) = 
R(A) + R(K(B)\setminus K(A))
\label{eq:1}
\end{equation}
We denote $K'(A) = \bigcup_{C\in K(A)}<C,b>$  and $ K'(D_{(n,k,k-1)}) =\bigcup_{C\in K(D_{(n,k,k-1)})}<C,b_1>,$ 
where   $b_1=(1,\ldots,1)^t$. By Corollary  \ref{con:veq} we have   
$
R(K(B)\setminus K(A)) = \frac{1}{n - k + 1}R(K'(A)). 
$
Substituting this   into  \eqref{eq:1}, we get 
\begin{equation} \label{eq:d1}
R(B) = R(A) + \frac{1}{n - k + 1}R(K'(A)).
\end{equation}
Repeating the same arguments for   $D_{(n,k,k-1)}$ and $D_{(n+1, k, k-1)}$, we obtain 
\begin{equation}  \label{eq:d2}
R(D_{(n+1, k, k-1)}) = R(D_{(n,k, k-1)}) + \frac{1}{n - k + 1}R(K'(D_{(n,k, k-1)})). 
\end{equation}

By assumptions $R(D_{(n,k,k-1)}) \lep R(A)$. Then we can construct the bijection between the matrices from $K(A)$ and  $K(D_{(n,k,k-1)})$ such that the matrix  $C \in K(A)$ corresponds to the matrix  $C' \in K(D_{(n,k,k-1)})$ satisfying  $\rk C \ge \rk C'$. Since any matrix  $C' \in K'(D_{(n,k,k-1)})$ has the form $D_{(k,\rk C' - 1)}$,  by Lemma  \ref{lem:n2} 
we have  the inequality  $R(<C',b_1>) \lep R(<C,b>)$.
Applying Lemmas \ref{lem:setsum} and \ref{lem:n1}, we obtain
\begin{equation}
 R(K'(D_{(n,k,k-1)})) = \suml{C' \in K(D_{(n,k,k-1)})}{}R(<C',b_1>) \lep 
\label{eq-n1}\end{equation}
$$
\lep \suml{C \in K(A)}{}R(<C,b>) = R(K'(A)).
$$

Using again that  $R(D_{(n,k,k-1)}) \lep R(A)$ we get the required result from Lemma  \ref{lem:n1} applied to the equalities \eqref{eq:d1} and \eqref{eq:d2}.

If $R(D_{(n,k,k-1)}) <^* R(A)$, then in the bijection constructed above we can find matrices $C \in K(A)$ and $C' \in K(D_{(n,k,k-1)})$ such that $\rk C' < \rk C$. It follows that $R(<C',b_1>) <^* R(<C,b>)$ and we obtain the strict inequality in \eqref{eq-n1} which implies $R(D_{(n+1,k,k-1)}) <^* R(B)$.
\end{prf}

Now we can prove the main result of this section.

\begin{lemma}
Let $A \in M_{k,n}(\pm 1)$ and $\rk A = k$. 
Then 
$R(D_{(n,k,k-1)}) \lep R(A)$. Moreover, if $R(D_{(n,k,k-1)}) = R(A)$, then there are $n-k+1$ columns in $A$ which are equal up to the multiplication by $(-1)$.
\label{lem:vec}\end{lemma}

\begin{prf}
Since  $\rk A=k$, there exists a nonsingular submatrix  $A_1$ of $A$ of the size $k$. Without loss of generality 
we can assume that it consists of the first  $k$ columns of $A$. Let us subsequently add the set of the remaining  columns of $A$
to  $A_1$ obtaining the matrices  $A_2, \ldots, A_{n-k+1}$, correspondingly. Note that  $A_{n-k+1}=A$. Similarly,  $D_{(n,k,k-1)}$ can be obtained 
from  $D_{(k,k-1)}$ by the subsequent addition of the columns $(1,\ldots,1)^t$. Each of the intermediate matrices is $D_{(k+i, k,k-1)}$, $i=1,\ldots, n-k$.  

For the nonsingular square matrices of the size  $k$ we have the equality
$
R(A_1) = (1, 0, \ldots, 0) = R(D_{(k,k-1)}).
$
Then  the conditions of Lemma  \ref{lem:vec0} are satisfied for the matrices  $A_1$ and $D_{(k,k-1)}$. 
Therefore, $R(D_{(k+1, k, k-1)}) \lep R(A_2)$. Applying subsequently  Lemma  \ref{lem:vec0} $(n-k)$ times we obtain the required inequality.

In order to prove the second part of the lemma we suppose that $R(D_{(n,k,k-1)}) = R(A)$. Let us prove by the induction on $n$. We can represent $A$ as the concatenation of the matrix $A' \in M_k(\pm 1)$ with $\rk A' = k$ and the remaining  $n-k$ columns. Induction step is to join one more column. 
The base of the induction is $n = k$, for which the statement is trivial. 

By induction hypothesis we assume that $R(D_{(n-1,k,k-1)}) = R(A'')$ for some $A''$ such that $A$ can be obtained from $A''$ by a concatenation with a column $b''$ and there are $n-k$ columns of $A''$  proportional to a certain vector $b'$. If $R(D_{(n-1,k,k-1)}) <^* R(A'')$, then by Lemma \ref{lem:vec0}  $R(D_{(n,k,k-1)}) <^* R(A)$, which is a contradiction. Therefore the equality in the formula \eqref{eq-n1} should be satisfied. By Corollary \ref{con1} the minimum of each summand in the left hand side  is achieved. However, in the right hand side the minimum can be achieved only if each matrix $C$ in the sum contains a column $c$ which is proportional to $b$. Since $n-k$ columns of $A''$ are proportional to $b'$, it follows that $b''$ is  proportional to $b'$ and there are $n-k+1$ proportional columns in $A$. The lemma is proved.  
\end{prf}

\begin{cons}
Let $A \in M_{k,n}(\pm 1)$ and $\rk A = k$, $R(A) = (r_i)$ and $R(D_{(n,k,k-1)}) = (u_i)$. Let $u_m$ be the last positive element in vector $(u_i)$. Then $r_m \le u_m$ 
and the equality is possible if and only if $R(A) = R(D_{(n,k,k-1)})$.
\end{cons}

\begin{prf}
The inequality $r_m \le u_m$ is obtained in the same way as in the proof of Corollary \ref{new-cons-2}. 

Suppose that $r_m = u_m$. Repeating arguments of the previous lemma we obtain that the equality in the formula \eqref{eq-n1} should be satisfied  for matrices of the rank $k-m+1$ and $k-m+2$. By Corollary \ref{con1} it is possible only if  there are $n-k+1$ proportional columns in the matrix $A$ and $R(A) = R(D_{(n,k,k-1)})$.
\end{prf}

We need the following special case of the last corollary.

\begin{cons}
Let $A \in M_{k,n}(\pm 1)$ and $\rk A = k$, $R(A) = (r_i)$ and $R(D_{(n,k,k-1)}) = (u_i)$. Then $r_k \le u_k$ 
and the equality is possible if and only if $R(A) = R(D_{(n,k,k-1)})$ or $r_k = u_k = 0$ and $n < 2k - 1$.
\label{new-cons}\end{cons}




\section{Proof of Kr\"auter conjecture for singular matrices} \label{S7}

In this section we use the  introduced order $\lep$ on  rank vectors to estimate the values of the functions $\per$ and $\mper$ (see Definition~\ref{notemper}). 
We obtain upper bounds for $\mper$ of  the full rank matrix $A$   via $\mper D_{(n,k,k-1)}$ for rectangular matrices $A$ with $k$ rows. In Lemma  \ref{lem:pr} the case $k>4$ is considered. The case $k=2$ is considered in Lemma  \ref{lem-pr-2}, and case $k=3$ is considered in Lemma \ref{lem-pr-3}. We investigate $k=4$ in  Lemma \ref{lem:28}  for $n=5$  and in Lemma \ref{lem:pr2} for $n>5$. Theorem \ref{max} combines the results of these 5 lemmas in the inequality \eqref{eq:centr3}, which is central for this section. Then we use inequality \eqref{eq:centr3} to prove   Conjecture \ref{h1} for singular matrices in Theorem~\ref{theorem2}.

\begin{note} By Corollary \ref{Cor3.3} and Theorem \ref{lem_n4} we get that Conjecture \ref{h1} holds true for all square matrices of size 2, 3, 4. \end{note}

\begin{lemma}
Let $A \in M_{k,n}(\pm 1)$,  
$\rk A = k>4$ and $n > k$. Assume that Conjecture \ref{h1} is proved for all square matrices of all sizes  $m \le k$. Then the following inequality holds:
\begin{equation}
\mper A \le \mper D_{(n, k, k-1)}.
\label{eq:centr}
\end{equation}
Moreover, the equality holds if and only if $A$ can be reduced to  $D_{(n,k,k-1)}$ by standard transformations.
\label{lem:pr}\end{lemma}

\begin{prf}
Since Conjecture \ref{h1} holds for $m \le k$ for any $C \in K(A)$,    $\mper C \le \mper D_{(k, \rk C - 1)}$. The equality is possible if and only if   $C$ can be reduced to $D_{(k, \rk C - 1)}$ by standard transformations. Let us sum such inequalities over all matrices $C \in K(A)$ and group summands depending on rank of $C \in K(A)$. We obtain the following inequality: \qquad $\mper A = $
$$
=\suml{C \in K(A)}{}\mper C = 
\suml{i = 1}{k}\suml{\begin{smallmatrix} C\in K(A): \\ \rk C = k-i+1 \end{smallmatrix}}{} \mper C \le \suml{i = 1}{k}\suml{\begin{smallmatrix} C\in K(A): \\ \rk C = k-i+1 \end{smallmatrix}}{} \mper D_{(k,k-i)}.
$$
Let $R(A) = (r_1,\ldots,r_k)$ and $R(D_{(n,k,k-1)}) = (r_1',\ldots,r_k')$. By Lemma \ref{lem:vec} $R(D_{(n,k,k-1)}) \lep R(A)$. Then we can construct a bijection $\rho:  K(A) \to K(D_{(n,k,k-1)})$ such than for each $C \in K(A)$   the  inequality $\rk C \ge \rk \rho(C)$ holds. By Conjecture \ref{h1} $\mper C \le \mper D_{(k, \rk C - 1)}$ and by Lemma \ref{lem40} $\mper D_{(k,\rk C - 1)} \le \mper D_{(k,\rk \rho (C) - 1)}$. Thus we obtain 
$$\mper C \le \mper D_{(k,\rk C - 1)} \le \mper D_{(k,\rk \rho (C) - 1)} = \mper \rho(C)$$
and the equality is obtained if and only if $C$ can be reduced to $\rho(C) = D_{(k,\rk \rho (C) - 1)}$ by standard transformations and $\rk C = \rk \rho(C)$. 

Now let us consider the difference between $\mper A$ and $\mper D_{(n,k,k-1)}$. Using the constructed bijection $\rho$ and the last inequality we obtain:
$$
\mper D_{(n,k,k-1)} - \mper A = \suml{C\in K(A)}{} (\mper \rho(C) - \mper C) \ge 0. 
$$

We proved the required inequality. Now we show that if the equality is obtained, then the matrix $A$ can be reduced to $D_{(n,k,k-1)}$ by standard transformations. 

The equality is obtained if and only if $C$ can be reduced to $\rho(C)$ by standard transformations. It follows that $\rk C = \rk \rho(C)$ and $\mper C = \mper \rho (C) = \mper D_{(k, \rk C - 1)}$ for every $C \in K(A)$. Thus we obtain $R(A) = R(D_{(n,k,k-1)})$.

By Lemma \ref{lem:vec} there exist  $n-k+1$ proportional columns in  $A$. Without loss of generality these are the columns with indices from  $k$th till  $n$th. 
Thus the submatrix $C' = A[1,\ldots, k| 1, \ldots, k]$ has rank $k$. Since $\rk C' = \rk \rho(C') = k$ and $\mper C' = \mper \rho (C') = \mper D_{(k,k-1)}$, and Conjecture \ref{h1} holds true for the matrices of size $k$ the inequality is possible only if $A[1,\ldots, k| 1, \ldots, k]$ can be reduced to  $D_{(k,k-1)}$ by standard transformations. 

Thus $A[1,\ldots, k| 1, \ldots, k] = D_{(k,k-1)}$. Since $k$th column of $A$ consists of the entries $1$ only, and the columns from   $k$th till $n$th are proportional, it is possible to multiply them by $-1$ to obtain  $A = D_{(n,k,k-1)}$, which concludes the proof. 
\end{prf}

Our next goal is to prove the analogs of Lemma \ref{lem:pr} for rectangular matrices with  $4, \; 3$, or $2$  rows. 

\begin{lemma}
Let $A \in M_{2,n}(\pm 1)$, $\rk A = 2 < n$. Then $\mper A \le \per D_{(n,2,1)}$. The equality holds if and only if $A$ can be reduced to   $D_{(n,2,1)}$ by standard transformations. 
\label{lem-pr-2}\end{lemma}

\begin{prf}
Let us multiply the columns of $A$ by  $(-1)$ in such a way that all entries in the second row are $1$.  Up to standard transformations without loss of generality we can denote this new matrix again by $A$. Let the first row of $A$ contain  $k$ negative entries. Since $\rk A = 2$, its rows are not proportional. Hence,  $1 \le k \le n-1$. Therefore by definition  $\mper A$ is a sum of $\frac12 n(n-1)$ summands such that   $k(n-k)$ summands are $0$, and the other summands are equal to $2$. Then the value of   $\mper A$ is maximal if and only if $k(n-k)$ is minimal, i.e., if and only if either   $k = 1$ or   $k = n-1$. In both of these cases it is straightforward to check that   $A$ can be reduced to   $D_{(n,2,1)}$ by standard transformations.  
\end{prf}

\begin{lemma}
Let $A \in M_{3,n}(\pm 1)$, $\rk A = 3$ and $n \ge 4$. Then $\mper A \le \mper D_{(n,3,2)}$. For $n \ge 5$ the equality holds if and only if $A$ can be reduced to $D_{(n,3,2)}$ by standard transformations. For $n = 4$ the equality holds if and only if $A$ can be reduced to $D_{(4,3,2)}$ or $D_{(4,3,3)}$ by standard transformations.
\label{lem-pr-3}\end{lemma}

\begin{prf}
Let $R(A) = (r_1, r_2, r_3)$. By Lemma   \ref{lem_n3} a $3\times 3$ matrix  $B$ can be reduced to   $D_{(3,\rk B - 1)}$ by standard transformations. Hence, 
$\mper A = \suml{B \in K(A)}{} |\per B| = 2 (r_1 + r_2) + 6r_3.$
Similarly, if $R(D_{(n,3,2)}) = (u_1, u_2, u_3)$, then we get
$\mper D_{(n,3,2)} = \suml{B \in K(D_{(n,3,2)})}{} |\per B| = 2 (u_1 + u_2) + 6u_3.$

By the conditions it holds that $r_1 + r_2 + r_3 = u_1 + u_2 + u_3$. Denote $S = r_1 + r_2 + r_3$. 
Then
$
\mper A  - \mper D_{(n,3,2)} = 2 (S - r_3) + 6r_3 - 2(S - u_3) - 6u_3 = 4 (r_3 - u_3).
$
By Lemma \ref{lem:vec} $R(A) \leq R(D_{(n,3,2)})$. Therefore $r_3 - u_3 = (S - r_1 - r_2) - (S - u_1 - u_2) = (u_1 - u_2) - (r_1 - r_2) \le 0$. Thus we proved the required inequality. 
 
Let $n \ge 5$. If $\mper A$ is maximal, then   $r_3 = u_3$. By Corollary \ref{new-cons} $A$ has $n-2$ proportional columns. Up to the standard transformations we may assume that all proportional columns are positive and are located on the positions from  $3$ until $n$. If the first or the second columns contain more than $1$ negative entry, then we multiply the corresponding column by $(-1)$. Note that we cannot obtain a positive column since in this case $\rk A \le 2$. Thus,  $A$ contains two negative entries lying in different columns. From the condition   $\rk A =3$ it follows that these entries lie in different rows also. Hence $A$ can be reduced to   $D_{(n,3,2)}$ by standard transformations.  

Let $n = 4$. Let us multiply columns of $A$ to obtain at most one negative element in each column. Denote the number of columns with negative elements by $l$. There are several options.

{\bf 1.}   $l = 0$ or $l = 1$. Then $\rk A < 3$. This contradicts the conditions. 

{\bf 2.}   $l = 2$. If $-1$s lie in the same row, then $\rk A < 3$. This contradicts the conditions. If $-1$s lie in different rows, then $A$ can be reduced to $D_{(4,3,2)}$.

{\bf 3.}   $l = 3$. If all $-1$s lie in different rows, then $A$ can be reduced  to $D_{(4,3,3)}$. 

If there are two $-1$s in the same row, then we multiply this row by $-1$. In the obtained matrix there are three $-1$s and two of them lie in the same column. Multiplying this column by $-1$ we obtain Case~{\bf 2}.

If all negative elements lie in the same row we multiply this row by $-1$. In the   obtained matrix there is one negative element and we are in Case~{\bf 1}.

{\bf 4.}  $l = 4$. If all $-1$s lie in the same row, then $\rk A < 3$.  This is a contradiction.

If three $-1$s lie in the same row, then we multiply this row by $-1$. In obtained matrix there are two $-1$s, and they lie in the same column. We can multiply this column by $-1$ and obtain Case~{\bf 1}.

Suppose at most two negative entries lie in the same row. We multiply this row by $-1$. In the obtained matrix there are four  $-1$s which lie in two columns, two $-1$s in each column. Let us multiply these columns by $-1$ to obtain Case {\bf 2} in which the lemma is proved already. 
\end{prf}

Our next goal is to prove the inequality  \eqref{eq:centr} for matrices from $M_{4,n}$ where $n > 4$. In order 
to show this, let us prove the following lemma. 
\begin{lemma}
Let $A \in M_{4,5}(\pm 1)$ and $\rk A = 4$. Then it holds that
\begin{equation}
\mper A = 
\suml{C \in K(A)}{} |\per C| \le \suml{C' \in K(D_{(5,4,3)})}{} |\per C'| = \mper D_{(5,4,3)}
\label{eq28}
\end{equation}
 The equality in \eqref{eq28} holds if and only if $A$ can be reduced to   $D_{(5,4,3)}$ by standard transformations.   \label{lem:28}\end{lemma}

\begin{prf}
Since $\rk A = 4$ and $\mper A$ does not depend on the order of columns of $A$, without loss of generality $\rk A(|5) = 4$.

Recall that by Theorem  \ref{experiment} if
$ X $ is a nonsingular matrix of size $4$ satisfying  $\per D_{(4,3)} < \per X$, then $X$ can be converted to $D_{(4,4)}$ by the standard transformations.
There are two subcases: 

{\bf 1.} Suppose $A$ does not contain $D_{(4,4)}$ as a submatrix up to standard transformations. Then for any submatrix $A_i$ of rank $k$ we obtain $\per A_i \le \per D_{(4,k-1)}$. It means that we can repeat the arguments from Lemma \ref{lem:pr}. It follows that $\mper A \le \mper D_{(5,4,3)}$ and the equality holds if and only if $A$ can be reduced to $D_{(5,4,3)}$ by standard transformations. 

{\bf 2.} Assume that some submatrix of $A$ coincides with the matrix  $D_{(4,4)}$ up to the standard transformations. Since $\rk D_{(4,4)} = 4$, without loss of generality $A(|5) = D_{(4,4)}$. 

Then $A$ has the   form $\left(\begin{smallmatrix}
-1 & 1 & 1 & 1 & x_1 \\
1 & -1 & 1 & 1 & x_2 \\
1 & 1 & -1 & 1 & x_3 \\
1 & 1 & 1 & -1 & x_4
\end{smallmatrix}\right)
,$ 
where $x_i \in \{\pm 1\}$. 

If we multiply the fifth column by $-1$, then  by definition $\mper A$ does not change. It means that up to the standard transformations we may assert that at least two elements of the last column are positive. 

Observe that if we apply the same permutation to the rows  $1, \; 2,\; 3,\; 4$ and columns  $1, \; 2,\; 3,\; 4$ of the matrix  $A$, then the block  $A(|5)$ remains invariant and the value  $\mper A$ remains unchanged. Therefore  we can apply such permutation in order to get $x_1 \le x_2 \le x_3 \le x_4$. Since  $x_i \in \{\pm 1\}$,  there are three different possibilities for the matrix  $A$
determined by the last column: 
\begin{enumerate}
\item All elements are positive, $A=A_1=\left(\begin{smallmatrix}
-1 & 1 & 1 & 1 & 1 \\
	1 & -1 & 1 & 1 & 1 \\
	1 & 1 & -1 & 1 & 1 \\
	1 & 1 & 1 & -1 & 1
	\end{smallmatrix}\right)$.
	
\item $x_1 = -1$ and $x_2 = x_3 = x_4 = 1$, $A=A_2=\left(\begin{smallmatrix}
	-1 & 1 & 1 & 1 & -1 \\
	1 & -1 & 1 & 1 & 1 \\
	1 & 1 & -1 & 1 & 1 \\
	1 & 1 & 1 & -1 & 1
	\end{smallmatrix}\right)$.
	
	\item $x_1 = x_2 = -1$ and $x_3 = x_4 = 1$, $A=A_3=\left(\begin{smallmatrix}
-1 & 1 & 1 & 1 & -1 \\
1 & -1 & 1 & 1 & -1 \\
1 & 1 & -1 & 1 & 1 \\
1 & 1 & 1 & -1 & 1
	\end{smallmatrix}\right)$.
\end{enumerate}
In each of these cases direct computations give
$
\mper (A_1)= 24$,
$ \mper (A_2)= 8 ,$ and $\mper (A_3)= 24.$ 
Since $\mper D_{(5,4,3)} = 32$, the lemma is proved.
\end{prf}

\begin{lemma}
Let  $A \in M_{4,n}(\pm 1)$,  $\rk A=4$, $R(A) = (r_i)$, and $n > 4$.
Let the statement of Conjecture \ref{h1} holds for the singular matrices of the order 4. 
Then the   inequality 
\begin{equation}
\mper A \le  
\mper D_{(n,4,3)}.
\label{eq:centr2}
\end{equation}
is true. The equality holds if and only if $A$ can be reduced to   $D_{(n,4,3)}$ by standard transformations. 
\label{lem:pr2}\end{lemma}

\begin{prf}
If there is no submatrix of $A$ which can be reduced to  $D_{(4,4)}$ by standard transformations, then the arguments from Lemma  \ref{lem:pr} work. Note that the maximal value of the $\mper A$ can be achieved in this case on the matrix equivalent to $D_{(n,4,3)}$ up to standard transformations. If there is such a submatrix, then we can without loss of generality assume that $A[1,2,3,4|1,2,3,4] = D_{(4,4)}$. 

By Lemma  \ref{lem:vec} the vector $R(D_{(n,4,3)}) \lep R(A)$. Then there exists a bijection $\rho: K(A)\to K(D_{(n,4,3)})$  such that  each $C \in K(D_{(n,4,3)})$ is the image of   $C' \in K(A)$ for which $\rk C'\ge \rk C$. Note that there is a unique matrix of rank 4 in  $K(D_{(n,4,3)})$, and this is  $D_{(4,3)}$. Without loss of generality we may assume that $\rho(A[1,2,3,4|1,2,3,4])=D_{(4,3)}$. By assumption  Conjecture  \ref{h1} holds for singular matrices of size $4$. Then using Proposition \ref{PropSmallD} we obtain that $\mper \rho(C') \ge \mper C'$ for all   $C' \in K(A)\setminus \{A[1,2,3,4|1,2,3,4]\}$. 

Let us consider the submatrix  $A[1,2,3,4|1,2,3,4,5]$. Its submatrices provide five elements of   $K(A)$. We set in correspondence to these elements one matrix of the rank $4$ and four matrices of the rank less than or equal to 3. Without loss of generality these are submatrices of $D_{(n,4,3)}[1,2,3,4|1,2,3,4,5] = D_{(5,4,3)}$. It is shown in the proof of  Lemma \ref{lem:28}   that $\mper A[1,2,3,4|1,2,3,4,5] \le 24$. By direct computations $\mper D_{(5,4,3)} = 32$. 

We split  $K(A) = K_A'\cup K_A''$ in such a way that $K_A'$ contains all submatrices of  $A[1,2,3,4|1,2,3,4,5]$. Similarly we split   $K(D_{(n,4,3)}) = K_D'\cup K_D''$, where  $K_D'$ contains all submatrices of  $D_{(n,4,3)}[1,2,3,4|1,2,3,4,5]$. By construction $\rho$ maps   $K_A'$ onto  $K_D'$ bijectively. Hence  there is a bijection $\rho''$  between   $K_A''$ and $K_D''$ satisfying the condition  $\rk \rho''(M) \le \rk M$ for all $M\in K_A''$ . Since Conjecture \ref{h1} holds for singular matrices of order  $4$  and for all  $C \in K_D''$ we have $\rk C \le 3$, it follows that Conjecture \ref{h1} holds for matrices from   $K_D''$. Hence by Proposition  \ref{PropSmallD} for   $M\in K_A''$   one of the following cases holds:

1.  $\rk M=4$. Then $\rk \rho (M) \le 3$. Thus $\mper M \le 8 \le \mper D_{(4,2)} \le \mper D_{(4,\rk \rho(M) -1)} = \mper \rho(M)$.

2. $\rk M\le3$. Then $\rk \rho (M) \le \rk M$ and by the induction hypothesis and Proposition \ref{PropSmallD} there are inequalities $\mper M \le \mper D_{(4,\rk M -1)} \le \mper D_{(4,\rk \rho(M) -1)} = \mper \rho(M)$.

Thus
$
\mper A=\suml{C'\in K(A)}{}\mper C' = \suml{C'\in K_A'}{}\mper C'  + \suml{C'\in K_A''}{}\mper C' = 24 + \suml{C'\in K_A''}{}\mper C' < 32 + \suml{C \in K_D''}{}\mper C = \suml{C \in K_D'}{}\mper C+ \suml{C \in K_D''}{}\mper C = \suml{C \in K(D_{(n,4,3)})}{}\mper C = \mper D_{(n,4,3)}
$.
The lemma is proved.
\end{prf}

Summarizing the results of Lemmas \ref{lem-pr-2}, \ref{lem-pr-3}, \ref{lem:pr2}, and  \ref{lem:pr} we get the following theorem:
\begin{theorem}
Let $A \in M_{k,n}(\pm 1)$, $k < n$,  $\rk A = k \ge 2$, and Conjecture \ref{h1} holds for matrices 
of the size  $l$ where $l \le k$. Then we have the following inequality:
\begin{equation}
\mper A \le \mper D_{(n,k,k-1)}.
\label{eq:centr3}
\end{equation}
For $k \ge 2$ and $(k,n) \ne (3,4)$ the inequality \eqref{eq:centr3} becomes equality if and only if $A$ can be reduced to $D_{(n,k,k-1)}$ by standard transformations. For $k = 3$ and $n = 4$ the inequality \eqref{eq:centr3} becomes equality if and only if $A$ can be reduced to $D_{(4,3,2)}$ or $D_{(4,3,3)}$.
\label{max}\end{theorem}

Now we are ready to prove Kr\"auter conjecture for singular matrices.

\begin{theorem}
Let  $A \in M_n(\pm 1)$,  $\rk A =k < n$. Assume that Conjecture  \ref{h1} is proved for all matrices of the size  $l$,
$l \le k$. Then we have the following inequality: 
$$
|\per A| \le \per D_{(n,k-1)}.
$$ 
Moreover, if $|\per A| = \per D_{(n,k-1)}$, then $A$ can be reduced to $D_{(n,k-1)}$ by standard transformations. 
\label{theorem2}\end{theorem}

\begin{prf}
By Theorem~\ref{th-kr-new} the statement holds if  $k = 1$. If $k = 2$ and $n = 3$ the result follows from Lemma~\ref{lem_n3}. 
Thus   we further assume   $k \ge 2$ and $n \ge 4$. 

Up to the standard transformations we consider the case: the first $k$ rows of $A$ are linearly independent. We use the Laplace decomposition formula for these rows. Let $\al = \{1, \ldots, k\}$. Then
$
|\per A| = \Bigl|\suml{\be \in \SS_{n,k}}{} \per A(\al|\be) \per A[\al|\be] \Bigr|   \le \suml{\be \in \SS_{n,k}}{} |\per A(\al|\be)| |\per A[\al|\be]|.
$
For any $\be \in \SS_{n,k}$ we have $\per A(\al|\be) \le (n-k)!$. By Theorem \ref{max} we get:
\begin{equation} \label{eq:bound} \begin{array}{c}
|\per A|  \le (n-k)! \suml{\be \in \SS_{n,k}}{} |\per A[\al|\be]| = (n-k)!\mper A[\al|1,\ldots, n] \le  
\\
\le  (n-k)! \cdot \mper D_{(n,k,k-1)} = |\per D_{(n,k-1)}|.
\end{array}
\end{equation}

Now we are going to characterize the equality case.

First we consider the exceptional case $k = 3$ and $n = 4$ of Theorem \ref{max} when the maximal value of $\mper A[\al|1,\ldots, n]$ can be obtained on $D_{4,3,3}$. Up to the standard transformations  $A$ is of the following form:
$
\left( \begin{smallmatrix}
x_1 & x_2 & x_3 & x_4 \\
-1 & 1 & 1 & 1 \\
1 & -1 & 1 & 1 \\
1 & 1 & -1 & 1 
\end{smallmatrix} \right).
$

By Laplace decomposition in the first row, the maximal value of $|\per A| = 8$ can be achieved if and only if $(x_1 , x_2 , x_3 , x_4) = \pm(1, 1, 1, -1)$. Then $A = D_{(4,4)}$ up to the standard transformations. Since $\rk D_{(4,4)} = 4$, this contradicts the condition $\rk A =k= 3$. Hence $|\per A| \le 4$ by Lemma~\ref{del}. 

Now we consider the general case. By Theorem \ref{max} $\mper A[\al|1,\ldots, n] = \mper D_{(n,k,k-1)}$ if and only if $A[\al|1,\ldots, n]$ can be reduced to $D_{(n,k,k-1)}$ by standard transformations. Thus we consider $A[\al|1,\ldots, n] = D_{(n,k,k-1)}$.

For any $\be \in \SS_{n,k}$ we have $|\per A(\al|\be)| \le (n-k)! = \per D_{(n,k,k-1)} = \per J_{n-k}$. 
If there is $\be' \in \SS_{n,k}$ such that $|\per A(\al|\be')| < (n-k)!$ it follows that in the first part of \eqref{eq:bound} there is strict inequality and $\mper A < \mper D_{(n,k-1)}$. Thus equality in \eqref{eq:bound} is possible only if $\rk A(\al|\be) = 1$ for any $\be \in \SS_{n,k}$. It follows that $\rk A(\al|) = 1$ since rank of any submatrix of maximal dimension is equal to $1$. 

For each $i > k$ if $a_{i1} = -1$ we multiply the $i$th row by $-1$. Since $\rk A(\al|) = 1$, all columns of $A(\al|)$ are equal to $\pm b$ where $b = (1, \ldots, 1)^t \in M_{n-k,1}(\pm 1)$.  
If all columns of $A(\al|)$ are equal to $b$, then $A = D_{(n,k-1)}$, and the theorem is proved. Suppose there is a column equal to $-b$ in $A(\al|)$. Then there are $\be_1, \be_2 \in \SS_{n,k}$ such that $\per A(\al|\be_1) = (n-k)!$ and $\per A(\al|\be_2) = -(n-k)!$. There are two options. 

{\bf 1.} By Lemma \ref{lem40} and Proposition \ref{PropSmallD} for any $k \ge 2$, $k > l \ge 0$ and $(k,l) \ne (2,1)$ we have the inequality $\per D_{(k,l)} > 0$. Thus in Laplace decomposition formula in \eqref{eq:bound} there are both positive and negative summands. Hence, the first inequality is strict and $|\per A| < \per D_{(n,k,k-1)}$.

{\bf 2.} Suppose $(k,l) = (2,1)$. If columns of   $A(\al|1)$ are equal to $-b$, then we multiply all rows of $A(\al|)$ and the first column of $A$ by $-1$. The obtained matrix contains only one negative element $a_{21} = -1$. Then it is permutationally equivalent to $D_{(n,2,1)}$ and the theorem is proved. If there are both positive and negative columns in $A(\al|1)$, then we can find the columns $i_1$ and $i_2$ in $A$ such that 
$A[1,2,3|1,i_1, i_2] = 
\left( \begin{smallmatrix}
-1 & 1 & 1 \\
1 & 1 & 1 \\
1 & 1 & -1
\end{smallmatrix} \right)=A'.
$
Since $\rk A' = 3$ we get $\rk A \ge 3$. This contradicts  $\rk A = k = 2$. 
\end{prf}

\begin{cons}
 Conjecture \ref{h1} is true.
\end{cons}

\begin{prf} We prove this conjecture by induction. The base of induction is given by Corollary \ref{Cor3.3} and  Theorem \ref{experiment}. The step of induction for non-singular matrices is proved  in Theorem \ref{theorem1}.   The step of induction for singular matrices is proved  in Theorem \ref{theorem2}. 
\end{prf}

Thus  Conjecture \ref{h1/1} is proved as well.

\section*{Acknowledgments}
The authors are grateful to the referees for valuable comments and suggestions.

\end{document}